\documentclass[11pt]{amsart}
\usepackage{amssymb}
\usepackage{hyperref}
\usepackage{verbatim}

\theoremstyle{plain}
\newtheorem{theorem}{Theorem}[section]

\newtheorem{proposition}[theorem]{Proposition}

\newtheorem{lemma}[theorem]{Lemma}

\theoremstyle{definition}

\newtheorem{remark}[theorem]{Remark}

\newtheorem{conjecture}[theorem]{Conjecture}
\newtheorem{problem}[theorem]{Problem}

\theoremstyle{remark}

\numberwithin{equation}{section}

\newcommand{\N}{\mathbb N}
\newcommand{\Z}{\mathbb Z}

\newcommand{\R}{\mathbb R}
\newcommand{\C}{\mathbb C}

\DeclareMathOperator{\GL}{GL}

\DeclareMathOperator{\Ot}{O}
\DeclareMathOperator{\SO}{SO}
\DeclareMathOperator{\SU}{SU}
\DeclareMathOperator{\Ut}{U}

\DeclareMathOperator{\Sp}{Sp}

\DeclareMathOperator{\Ad}{Ad}

\DeclareMathOperator{\spec}{Spec}
\DeclareMathOperator{\Hom}{Hom}

\DeclareMathOperator{\diag}{diag}

\DeclareMathOperator{\Ehr}{Ehr}

\newcommand{\op}{\operatorname}

\newcommand{\ba}{\backslash}
\newcommand{\mi}{\textrm{i}}
\newcommand{\norma}[1]{\|{#1}\|_1}

\newcommand{\HH}{\mathcal H}
\newcommand{\PP}{\mathcal P}

\newcommand{\zz}{\ell}

\title[Spectra of lens spaces]{Recent results on the spectra of lens spaces}

\author{Emilio A. Lauret}
\address{INMABB (CONICET) and Dpto.\ de Matemática, Universidad Nacional del Sur, Bah\'ia Blanca, Argentina.}
\email{emilio.lauret@uns.edu.ar}

\author{Roberto J. Miatello}
\address{CIEM--FaMAF (CONICET), Universidad Nacional de C\'ordoba, Medina Allende s/n, Ciudad Universitaria, 5000 C\'ordoba, Argentina.}
\email{miatello@famaf.unc.edu.ar}

\author{Juan Pablo Rossetti}
\address{CIEM--FaMAF (CONICET), Universidad Nacional de C\'ordoba, Medina Allende s/n, Ciudad Universitaria, 5000 C\'ordoba, Argentina.}
\email{rossetti@famaf.unc.edu.ar}

\thanks{This research was supported by grants from CONICET, FONCyT and SeCyT. The first named author was also supported by the Alexander von Humboldt Foundation (return fellowship)}
\date{September 2019}

\begin{document}

\begin{abstract}
In this paper we report on recent results by several authors, on the spectral theory of lens spaces and orbifolds and similar locally symmetric spaces of rank one.
Most of these results are related to those obtained by the authors in [IMRN (2016), 1054--1089], where the spectra of lens spaces were described in terms of the one-norm spectrum of a naturally associated congruence lattice.
As a consequence, the first examples of Riemannian manifolds isospectral on $p$-forms for all $p$ but not strongly isospectral were constructed.

We also give a new elementary proof in the case of the spectrum on functions.
In this proof, representation theory of compact Lie groups is avoided and replaced by the use of Molien's formula and a manipulation of the one-norm generating function associated to a congruence lattice.
In the last four sections we present several recent results, open problems and conjectures on the subject.
\end{abstract}

\maketitle

\tableofcontents

\section{Introduction}
This article surveys recent results on the spectral geometry of lens and related spaces.
In \cite{LMR-survey}, the authors discussed the seminal work of Ikeda (\cite{IkedaYamamoto79,Ikeda80_isosp-lens, Ikeda88}), together with results in \cite{LMR-onenorm}.
A number of new results have appeared since then, and a purpose of the present article is to give a rather complete update.

In this paper we study the spectrum of a canonical operator $\Delta_{\tau,\Gamma}$ acting on sections of a natural bundle on a space of the form $\Gamma\ba G/K$, where $G$ is a compact Lie group, $K$ is a closed subgroup of $G$, $\tau$ is a finite dimensional representation of $K$ and $\Gamma$ is a finite subgroup of $G$.
A main case is when $\Gamma\ba G/K$ is a lens space or a lens orbifold.

In one result, we give an elementary proof of the theorem in \cite{LMR-onenorm} concerning the spectrum of the Laplace--Beltrami operator acting on functions on lens spaces.
In particular, this approach avoids the use of irreducible representations of compact Lie groups (see Section~\ref{sec:spec0lentes} and also \cite{MohadesHonari17}).
An important tool is an integral lattice of rank $n$, defined by a congruence relation, naturally associated to each $(2n-1)$-dimensional lens space.

The spectrum of the Hodge--Laplace operator acting on smooth $p$-forms of an odd-dimensional sphere is considered in Section~\ref{sec:p-spectrum}.
An emphasis is given to isospectral characterizations and to the construction of isospectral examples based on the construction of one-norm isospectral lattices.
Results and examples from \cite{LMR-onenorm}, \cite{DeFordDoyle14}, and  \cite{Lauret-p-spectralens} are discussed in quite some detail.

We review the description of the spectrum of a lens space by using the Ehrhart theory for counting integer points in a naturally associated polytope. Furthermore, we include  applications due to Mohades and Honari~\cite{MohadesHonari17}, on
the toric varieties associated to pairs of isospectral lens spaces (see Section~\ref{sec:Ehrhart}).

The last section contains brief discussions on several recent articles related to the spectral theory of lens spaces.
Namely, the work on the spectrum of the Dirac operator on spin lens spaces by Boldt~\cite{Boldt15} and  Boldt-Lauret~\cite{BoldtLauret-onenormDirac};
the computational studies on $p$-isospectral lens spaces (\cite{GornetMcGowan06,Lauret-computationalstudy}), %the work on individual $p$-spectra of lens spaces (\cite{Lauret-p-spectralens}), 
the work by Bari and Hunsicker (\cite{Shams11,ShamsHunsicker17}) on the spectra of lens orbifolds, the article \cite{Lauret-spec0cyclic} extending the one-norm method to other compact symmetric spaces of real rank one, and the harmonic counting measure introduced 
in~\cite{MohadesHonari16}.

\subsection*{Acknowledgments}  We dedicate this paper to Manfredo Do Carmo, who was an excellent mathematician and a prominent developer of differential geometry in South America.
This article expands the talk given by the second named author at the ICM-satellite conference in July of 2018 in S\~ao Paulo. He wishes to  thank the organizers C.~Gorodski and P.~Piccione, for putting up together a high level and most pleasant conference.
 The authors are greatly indebted to the referee for many accurate comments and clarifications that have helped to improve significantly the presentation of the paper.

\section{Spectra on functions of spherical space forms}\label{sec:spectraspherical}

In this section we describe, by elementary methods, the spectrum of the Laplace operator acting on functions of a quotient of the sphere by a finite group acting by isometries.

We will use along this section some well-known facts on the spectrum of the unit sphere $S^{n-1}$ in $\R^n$, which can be found in \cite[Ch.\ 2.G and Ch.\ 3.C]{BGM}.
For $f\in C^\infty(S^{n-1})$, let $\widehat f:\R^{n}\smallsetminus \{0\}\to \C$ be given by $\widehat f(x)= f(\frac{x}{|x|})$.
For $g\in C^\infty(\R^{n}\smallsetminus \{0\})$, set $\widetilde g=g|_{S^{n-1}}$.
It is well known (see for instance \cite[Thm.~22.1]{Shubin}) that if  $\Delta$ is the Laplace--Beltrami operator on $S^{n-1}$, then $\Delta f= (\Delta_{\R^n} \widehat f)|_{S^{n-1}}$ for every $f\in C^\infty(S^{n-1})$, and
$\Delta \widetilde g= {\widetilde{\Delta_{\R^n} g}}$ for every $g\in C^\infty(\R^{n}\smallsetminus \{0\})$  constant along the lines through the origin.
It is well known (see \cite[p.\ 159]{BGM}) that if $g$ is a homogeneous polynomial (in $n$ variables) of degree $k$ with complex coefficients which is harmonic, then
\begin{equation}
\Delta \widetilde g = k(k+n-2)\, \widetilde g.
\end{equation}
In other words, such a function is an eigenfunction of $\Delta$ with eigenvalue $k(k+n-2)$.
These functions are known as \emph{spherical harmonics}.
We shall see that there are no other eigenfunctions of $\Delta$ on $S^{n-1}$.

Let
\begin{equation}\label{eq:PP_kHH_k}
\PP_k =\left\{
\text{homogeneous complex polynomials of degree } k    \right\} \text{ and }\
\HH_k =\left\{g\in \PP_k: \Delta_{\R^n}g=0\right\}.
\end{equation}
It is also well known (see \cite[p.\ 160]{BGM})
that $\PP_k= \HH_k\oplus |r|^2 \PP_{k-2}$, where $|r|^2 = x_1^2+\dots+x_n^2$,
that is, every polynomial $g\in \PP_k$ can be written in a unique way as $g=g_0+|r|^2 g_1$ with $g_0\in \HH_k$ and $g_1\in \PP_{k-2}$.
In particular, the restriction to $S^{n-1}$ of every polynomial (not necessarily homogeneous) can be written as a sum of
restrictions to $S^{n-1}$ of elements in $\HH_j$, $j\in\N_0$.

By the Stone--Weierstrass theorem,
it follows that the subspace $\bigoplus_{k\geq0} {\widetilde{\HH_k}}$ is a dense subspace of $L^2(S^{n-1})$.
Hence, the spectrum of $\Delta$ is given by the collection of eigenvalues $\lambda_k:=k(k+n-2)$, for $k\geq 0$, each of them with
multiplicity $\dim {\widetilde{\HH_k}}$ ($=\dim \HH_k$). Moreover, it is easy to see that $\dim \PP_k = \binom{k+n-1}{n-1}$ and furthermore $\dim \HH_k=\dim \PP_k-\dim \PP_{k-2}=\binom{n+k-2}{k}\frac{2k+n-2}{k+n-2}$.

Clearly, for homogeneous polynomials $g_1,g_2$ of degree $k$ and in $n$ variables, $g_1=g_2$ if and only if $\widetilde g_1=\widetilde g_2$.
Thus, it will cause no confusion if we avoid the tilde symbol on $\PP_k$ and $\HH_k$ to identify functions on $S^{n-1}$.

We now study the spaces covered by $S^{n-1}$.
The isometry group of the unit sphere $S^{n-1}$ is the orthogonal group $\Ot(n)=\{g\in \GL(n,\R): g^tg=I_n\}$, where the action of an element $g\in\Ot(n)$ on $x\in S^{n-1}=\{(x_1,\dots,x_n)^t\in \R^n:x_1^2+\dots+x_n^2=1\}$	is given by multiplication on the left, that is, $gx$.
The group of orientation-preserving isometries is $\SO(n)=\{g\in \Ot(n): \det(g)=1\}$.

A \emph{spherical space form} is a space of the form $\Gamma\ba S^{n-1}$, where $\Gamma$ is a finite subgroup of $\Ot(n)$ acting freely on $S^{n-1}$.
In this case, $\Gamma$ is isomorphic to the fundamental group of the space  provided $n\geq 3$.
A comprehensive study and classification of spherical space forms is given in \cite{Wolf-book}.
Even dimensional spheres, $S^{n-1}$ with $n$ odd, yield only $S^{n-1}$ itself and  $P^{n-1}(\R)$, while
odd dimensional spheres yield infinitely many space forms.

If $\Gamma$ is allowed to be any finite subgroup of $\Ot(n)$, acting freely or not, then $\Gamma\ba S^{n-1}$ may be a manifold or just a \emph{good orbifold}
(also called \emph{global orbifold}).
See \cite{Gordon12-orbifold} for a detailed discussion on the spectral theory of the Laplace operator associated to an arbitrary orbifold.
In either case, the dimension of the eigenspace corresponding to the eigenvalue $\lambda$ in the quotient $\Gamma\ba S^{n-1}$ equals the dimension of the $\Gamma$-invariants in the eigenspace with eigenvalue $\lambda$ in $S^{n-1}$.
Thus, the determination of $\spec(\Gamma\ba S^{n-1})$ does not depend on whether $\Gamma\ba S^{n-1}$ is a manifold or not.

We point out that we will always consider $\Gamma$ inside $\SO(n)$ and hence, $\Gamma\ba S^{n-1}$ is orientable.

The group $G=\SO(n)$ acts on functions on $S^{n-1}$ by $(g\cdot f)(x) = f(g^{-1}x)$ for $g\in G$ and $f\in C^\infty(S^{n-1})$.
If $\Gamma$ is a finite subgroup of $G$, the functions on $\Gamma\ba S^{n-1}$ correspond with the
functions on $S^{n-1}$ that are invariant by $\Gamma$.
Consequently,
\begin{equation}\label{eq:specGammabaG}
\spec(\Gamma\ba S^{n-1}) = \big\{\!\big\{ \underbrace{\lambda_k,\dots,\lambda_k}_{\dim \HH_k^\Gamma \text{-times}} : k\geq0 \big\}\!\big\},
\end{equation}
where $\HH_k^\Gamma = \{f\in {\HH}_k: \gamma\cdot f=f\text{ for all } \gamma\in\Gamma\}$.
In other words, the multiplicity of $\lambda_k$ in $\spec(\Gamma\ba S^{n-1})$ equals $\dim {\HH}_k^\Gamma$, which, naturally, could be equal to zero.

Ikeda and Yamamoto~\cite{IkedaYamamoto79} introduced the \emph{spectral generating function} $F_\Gamma(z)$ associated to $\Gamma\ba S^{n-1}$ given by
\begin{equation}\label{eq:spectralgeneratingfunction}
F_\Gamma(z) = \sum_{k\geq0} \dim {\HH}_k^\Gamma \; z^k.
\end{equation}
This formal power series encodes the spectrum of the Laplace--Beltrami operator of $\Gamma\ba S^{n-1}$ in the sense that
if $\Gamma$ and $\Gamma'$ are two finite subgroups of $G$, then $\Gamma\ba S^{n-1}$ and $\Gamma'\ba S^{n-1}$ are isospectral (i.e.\ $\spec(\Gamma\ba S^{n-1})=\spec(\Gamma'\ba S^{n-1})$) if and only if $F_{\Gamma}(z)=F_{\Gamma'}(z)$.
Thus, the problem of describing the spectrum of the Laplace--Beltrami operator of $\Gamma\ba S^{n-1}$ is equivalent to the combinatorial problem of finding an explicit expression for $F_{\Gamma}(z)$.

Since $\PP_k \simeq   \HH_k\oplus \PP_{k-2}$, one has that $\dim {\HH}_k^\Gamma = \dim {\PP}_k^\Gamma - \dim {\PP}_{k-2}^\Gamma$, hence
\begin{align}
F_{\Gamma}(z)
&= \sum_{k\geq0} (\dim {\PP}_k^\Gamma - \dim {\PP}_{k-2}^\Gamma)\, z^k
= \sum_{k\geq0} \dim {\PP}_k^\Gamma\, z^k - z^2\sum_{k\geq0} \dim {\PP}_{k-2}^\Gamma\, z^{k-2}\\
&= (1-z^2) \sum_{k\geq0} \dim {\PP}_k^\Gamma \, z^k. \notag
\end{align}
Now, Molien's formula (\cite{Molien}) states that $\displaystyle \sum_{k\geq0} \dim {\PP}_k^\Gamma \, z^k = \frac 1{|\Gamma|} \sum_{\gamma\in \Gamma} \frac{1}{\det(I_n-z\gamma)}$,
hence
\begin{equation}\label{eq:Molien}
F_{\Gamma}(z) = \frac{1-z^2}{|\Gamma|} \sum_{\gamma\in \Gamma} \frac{1}{\det(I_n-z\gamma)},
\end{equation}
where $I_n$ denotes the $n\times n$ identity matrix.
Formula \eqref{eq:Molien} was given by Ikeda in \cite[Thm.~2.2]{Ikeda80_3-dimI}.
On the other hand, later, Wolf in \cite[(2.4)]{Wolf01} noticed that \eqref{eq:Molien} is a consequence of Molien's formula.

\section{Spectra on functions of lens spaces}\label{sec:spec0lentes}

The goal of this section is to give a simple proof, using formula \eqref{eq:Molien}, of Theorem 3.6(i) in \cite{LMR-onenorm}, which shows that two lens spaces are isospectral (on functions) if and only if their corresponding congruence lattices are isospectral with respect to the one-norm $\norma{\cdot}$.
Such result was proven in \cite{LMR-onenorm} by using representation theory of compact Lie groups, which will be avoided in this section.
An alternative proof was given by Mohades and Honari in \cite{MohadesHonari17}.

We will restrict our attention to odd-dimensional spheres, say $S^{2n-1}$ for $n\geq2$.
A lens space has cyclic fundamental group and is of the form $L(q;s_1,\dots,s_n) := \langle \gamma \rangle \ba S^{2n-1}$, where
\begin{equation}\label{eq:gamma}
\gamma=\begin{pmatrix}
\cos (\frac{2\pi s_1}{q}) & \sin (\frac{2\pi s_1}{q}) \\
-\sin (\frac{2\pi s_1}{q}) & \cos (\frac{2\pi s_1}{q}) \\
&&\ddots\\
&&&& \cos (\frac{2\pi s_n}{q}) & \sin (\frac{2\pi s_n}{q}) \\
&&&&-\sin (\frac{2\pi s_n}{q}) & \cos (\frac{2\pi s_n}{q})
\end{pmatrix}
\end{equation}
for some $q\in\N$ and $s_1,\dots,s_n\in \Z$ are such that $\gcd(q,s_j)=1$ for all $1\leq j\leq n$.
If we relax the last condition and assume that $\gcd(q,s_1,\dots,s_n)=1$ (to ensure that $\gamma$ has order $q$), then the action of $\langle \gamma\rangle$ on $S^{2n-1}$ is not necessarily free, and the quotient is called a \emph{lens orbifold}, or an \emph{orbifold lens space}.

Lens spaces play a fundamental role in topology.
For instance, they provided the first pair of non-homeomorphic homotopy equivalent topological spaces.
The \emph{Reidemeister torsion} was introduced in 1930 to distinguish such pairs topologically.
Later, it was conjectured in 1971 by Ray and Singer, and then proved in 1978 by Cheeger and M\"uller, that this topological invariant has an analytic analogue called the \emph{analytic torsion}.

By applying formula \eqref{eq:Molien}, the spectral generating function of a lens orbifold $L:=L(q;s_1,\dots,s_n)=\Gamma\ba S^{2n-1}$ is given by
\begin{equation}\label{eq:spectralgeneratingfunctionlens}
F_L(z):= F_\Gamma(z) = \frac{1-z^2}{q} \sum_{l=0}^{q-1} \frac{1}{\prod_{j=1}^n (1-\xi_q^{ls_j}z)(1-\xi_q^{-ls_j}z)},
\end{equation}
where $\xi_q=e^{2\pi \mi /q}$.
This formula first appeared in the literature in \cite[Thm.~3.2]{IkedaYamamoto79}.
Later, Ikeda used this expression to give families of isospectral lens spaces of dimension $\geq 5$ (\cite{Ikeda80_isosp-lens}).
Bari~\cite{Shams11} extended his construction to lens orbifolds.
In the opposite direction, Yamamoto~\cite{Yamamoto80} proved the non-existence of $3$-dimensional isospectral lens spaces.
Bari and Hunsicker~\cite{ShamsHunsicker17} have recently considered the same problem in the context of $3$-dimensional lens orbifolds (see Subsection~\ref{subse:BariHunsicker}).

Next, we will show a direct relation between the spectrum of the Laplace--Beltrami operator on a lens space
(via the spectral generating function defined in \eqref{eq:spectralgeneratingfunction})
and the one-norm spectrum of an associated lattice.

It is straightforward to see that
\begin{equation}
\frac{1-z^2}{(1-\xi_q^{ls}z)(1-\xi_q^{-ls}z)}
= \frac{1}{1-\xi_q^{ls}z} + \frac{1}{1-\xi_q^{-ls}z} - 1
= \sum_{m\in \Z} \xi_q^{lms} z^{|m|}
\end{equation}
for every $l,s\in\Z$.
Thus, in light of \eqref{eq:spectralgeneratingfunctionlens}, we have
\begin{align}
F_L(z)
&= \frac 1{{(1-z^2)}^{n-1}} \, \frac1q \, \sum_{l=0}^{q-1} \, \prod_{j=1}^n \frac{1-z^2}{(1-\xi_q^{ls_j}z)(1-\xi_q^{-ls_j}z)} \notag\\
&= \frac 1{{(1-z^2)}^{n-1}} \, \frac1q \,\, \sum_{l=0}^{q-1} \,\, \prod_{j=1}^n \,\, \sum_{a_j\in\Z}  \xi_q^{la_js_j}z^{|a_j|} \notag \\
&= \frac 1{{(1-z^2)}^{n-1}} \, \sum_{(a_1,\dots,a_n)\in\Z^n} \left( \frac1q \sum_{l=0}^{q-1}  \xi_q^{l (a_1s_1+\dots+a_ns_n)}\right) z^{|a_1|+\dots+|a_n|} \notag\\
&= \frac 1{{(1-z^2)}^{n-1}} \, \sum_{(a_1,\dots,a_n)\in\Z^n :\;  q | a_1s_1+\dots+a_ns_n }  z^{\norma{(a_1,\dots,a_n)}} \notag
\end{align}
where in the last identity we have used that for any $k\in\Z$
\begin{equation}
\sum_{l=0}^{q-1} \xi_q^{lk}
=\begin{cases}q &\text{ if $q$ divides $k$},\\
0 &\text{ otherwise}.
\end{cases}
\end{equation}
As usual, we keep the notation $\norma{(a_1,\dots,a_n)} = |a_1|+\dots+|a_n|$.

We observe that the points in $\Z^n$ contributing to the last sum in the long computation form a lattice, and this allows us to make a reformulation of the identity.
Thus, it is natural to associate to a lens orbifold $L(q;s_1,\dots,s_n)$ the \emph{congruence lattice}
\begin{equation}\label{eq:congruence-lattice}
\mathcal L(q;s_1,\dots,s_n):= \{(a_1,\dots,a_n)\in\Z^n: a_1s_1+\dots+a_ns_n\equiv0\pmod q\},
\end{equation}
as in \cite[Def.~3.2]{LMR-onenorm}.

For $k$ a non-negative integer, we denote by $N_{\mathcal L}(k)$ the number of elements in $\mathcal L :=\mathcal L(q;s_1,\dots,s_n)$ with one-norm equal to $k$, that is,
\begin{equation}
N_{\mathcal L}(k) =\#\{\mu=(a_1,\dots,a_n)\in\mathcal L: \norma{\mu}=|a_1|+\dots+|a_n|=k\}.
\end{equation}
We define the \emph{one-norm theta function} of $\mathcal L$ by
\begin{equation}\label{eq:one-normtheta}
\vartheta_{\mathcal L}(z) := \sum_{k\geq0}N_{\mathcal L}(k)\, z^k.
\end{equation}

Since
$\vartheta_{\mathcal L}(z)=\sum_{\mu\in\mathcal L} z^{\norma{\mu}}$,
we have thus proved

\begin{theorem}\label{thm:Fvartheta}
If $L$ is a lens orbifold and $\mathcal L$ its associated congruence lattice, then we have that
\begin{equation}\label{eq:F_L(z)}
F_{L}(z) = \frac{\vartheta_{\mathcal L}(z) }{(1-z^2)^{n-1}}.
\end{equation}
\end{theorem}

This expression allows us to give a formula for the multiplicity of the eigenvalue $\lambda_k$ in the spectrum of the Laplace--Beltrami operator on $L$ in terms of the numbers $N_{\mathcal L}(k)$ for $k\geq0$.
Indeed, by using the expansion $(1-z^2)^{-(n-1)} = \sum_{r\geq0} \binom{r+n-2}{n-2}z^{2r}$, Theorem~\ref{thm:Fvartheta} implies that
\begin{align}\label{eq:F_Gamma(z)pto}
F_{\Gamma}(z)
&= \left(\sum_{r\geq0} \binom{r+n-2}{n-2}z^{2r}\right) \left( \sum_{k\geq0}N_{\mathcal L}(k)\, z^k\right) \\
&= \sum_{k\geq0} z^k \; \sum_{r=0}^{\lfloor k/2\rfloor} \binom{r+n-2}{n-2} N_{\mathcal L}(k-2r). \notag
\end{align}
In other words, the multiplicity of $\lambda_k=k(k+2n-2)$ in $\spec(L)$ is given by
\begin{equation}\label{eq:dimHH_k^Gamma}
\dim \HH_k^\Gamma = \sum_{r=0}^{\lfloor k/2\rfloor} \binom{r+n-2}{n-2} N_{\mathcal L}(k-2r),
\end{equation}
which is the expression in \cite[(20)]{LMR-onenorm}.
This gives a simple formula for the multiplicities of the eigenvalues on a lens space, obtained by elementary methods.

In view of Theorem~\ref{thm:Fvartheta}, it is now straightforward to obtain the following equivalence condition for isospectrality of lens spaces.

\begin{theorem}\cite[Thm.~3.6(i)]{LMR-onenorm} \label{thm:3.6(i)}
Let $L$ and $L'$ be two lens spaces and let $\mathcal L$ and $\mathcal L'$ be their corresponding congruence lattices.
Then, $L$ and $L'$ are isospectral if and only if $\mathcal L$ and $\mathcal L'$ are $\norma{\cdot}$-isospectral (i.e.\ $\vartheta_{\mathcal L}(z) = \vartheta_{\mathcal L'}(z)$).
\end{theorem}

\section{Spectra on sections of homogeneous vector bundles} \label{sec:approach}

In this section we introduce a general approach based on representation theory of compact groups to study the spectra of compact locally homogeneous spaces.
It will provide us the adequate tools to work on the following more general settings
\begin{itemize}
\item arbitrary normal homogeneous spaces (not only odd-dimensional spheres);
\item Laplace type operators acting on sections of vector bundles (not only on functions).
\end{itemize}
We end this section by showing how to apply this approach to recover the results from Section~\ref{sec:spec0lentes}.
This method was used in \cite{LMR-onenorm} to study the spectra of the Hodge--Laplace operator acting on smooth $p$-forms on lens spaces.
These results will be discussed in Section~\ref{sec:p-spectrum}.

\subsection{Vector bundles on homogeneous spaces}\label{subsec:E_tau}
Let $G$ be a compact Lie group and $K$ a closed subgroup of $G$.
We consider on $G/K$ a $G$-invariant metric induced by an $\Ad(G)$-invariant inner product $\langle \cdot,\cdot\rangle$ on the Lie algebra $\mathfrak g$ (e.g.\ any negative multiple of the Killing form when $G$ is semisimple).
The corresponding homogeneous Riemannian manifold $G/K$ is called \emph{normal}.

Any complex $G$-homogeneous vector bundle on $G/K$ is constructed as follows  (see \cite[\S{}5.2]{Wallach-book}):
for $(\tau,W_\tau)$ a unitary representation of $K$, we set $E_\tau:=G\times_\tau W_\tau = G\times W_\tau/\sim$ where $(gk,w)\sim (g,\tau(k)w)$ for every $g\in G$, $w\in W_\tau$, $k\in K$.
Let the corresponding projection  $: E_\tau\to G/K$,  be given by $[g,w]\mapsto gK$, where $[g,w]$ denotes the class of $(g,w)$ in $E_\tau$.
By the Peter-Weyl theorem and Frobenius reciprocity, the space of $L^2$-sections decomposes as a $G$-module as
\begin{equation}\label{eq:PeterWeyl}
L^2(E_\tau)\simeq \bigoplus_{\pi\in\widehat G} V_\pi\otimes \Hom_K(V_\pi,W_\tau).
\end{equation}
Note that the Hilbert sum over $\pi$ is restricted to the set
\begin{equation}
\widehat G_\tau:= \{\pi\in \widehat G: \Hom_K(V_\pi,W_\tau)\neq0\}.
\end{equation}
 Here, to the element $v\otimes A\in V_\pi\otimes \Hom_K(V_\pi,W_\tau)$ one associates the smooth section $G/K\to E_\tau$ given by $xK\mapsto [x,A(\pi(x^{-1})\cdot v)]$.

We consider the  elliptic operator $\Delta_\tau$ acting on smooth sections of $E_\tau$ by the  negative of the Casimir element $C$ with respect to the inner product $\langle\cdot,\cdot\rangle$.  As in  \cite[\S{}5.6]{Wallach-book} we call it the Laplace operator on the bundle $E_\tau$.
Since $C$ lies in the center of the universal enveloping algebra, the operator $\pi(-C): V_\pi\to V_\pi$ commutes with $\pi(g)$ for any $g\in G$.
By Schur's Lemma, for $\pi\in \widehat G$, $\pi( -C)$ acts by a scalar $\lambda(C,\pi)$ on $V_\pi$.
Moreover, $\lambda(C,\pi)\geq0$ and it can be expressed in terms of the highest weight $\lambda_\pi$ of $\pi$.
Indeed,
\begin{equation}\label{eq:Casimireigenvalue}
\lambda(C,\pi)=   \langle\lambda_\pi +\rho, \lambda_\pi +\rho \rangle -\langle \rho,\rho \rangle
= \langle \lambda_\pi, \lambda_\pi+ 2\rho \rangle,
\end{equation}
where $\rho$ equals half the sum of the positive roots of $(\mathfrak g_\C, \mathfrak t_\C)$ and where $\mathfrak t$ is the Lie algebra of a maximal torus $T$ of $G$.
From \eqref{eq:PeterWeyl}, we conclude that the spectrum of $\Delta_\tau$ is given
 by the numbers of the form $\lambda(C,\pi)$ for some $\pi\in\widehat G_{\tau}$,
and furthermore, for such $\lambda$, its multiplicity is given by
\begin{equation}\label{eq:multiplicitylambda}
\sum_{\pi\in\widehat G_\tau:\; \lambda(C,\pi)=\lambda} \dim V_\pi\; \dim \Hom_K(V_\pi,W_\tau).
\end{equation}

Let $\Gamma$ be a finite subgroup of $G$.
Clearly, $\Gamma$ acts on $E_\tau$ by $\gamma\cdot [g,w]=[\gamma g,w]$, thus $\Gamma\ba E_\tau$ is a vector bundle over $\Gamma\ba G/K$.
Furthermore,
\begin{equation}\label{eq:PeterWeylGamma}
L^2(\Gamma\ba E_\tau)\simeq \bigoplus_{\pi\in\widehat G} V_\pi^\Gamma \otimes \Hom_K(V_\pi,W_\tau)
\end{equation}
as complex vector spaces.
Consequently, the spectrum of the induced operator $\Delta_{\tau,\Gamma}$ acting on sections of $\Gamma\ba E_\tau$ can be described as follows.
\begin{proposition} \label{prop:tau-spectrum}
Let $G$ be a compact Lie group, $K$ a closed subgroup of $G$, $\Gamma$ a finite subgroup of $G$ and $\tau$ a finite dimensional representation of $K$.
Every eigenvalue $\lambda$ in the spectrum of $\Delta_{\tau,\Gamma}$ is of the form $\lambda(C,\pi)$ for some $\pi\in\widehat G_{\tau}$, with multiplicity given by
\begin{equation}\label{eq:multeigenvalue}
\sum_{\pi\in\widehat G_\tau:\; \lambda(C,\pi)=\lambda} \dim V_\pi^\Gamma\; \dim \Hom_K(V_\pi,W_\tau).
\end{equation}
\end{proposition}

Summing up, one needs the following two ingredients to determine the spectrum of the operator $\Delta_{\tau,\Gamma}$
\begin{itemize}
\item $\dim \Hom_K(V_\pi,W_\tau)$ for each $\pi\in \widehat G$ (which in particular determines  $\widehat G_\tau$)

\item the dimension $\dim V_\pi^\Gamma$ of the space of $\Gamma$-invariants in $V_\pi$,  for each $\pi\in\widehat G_\tau$.
\end{itemize}

The first ingredient involves the \emph{branching law} from $G$ to $K$, which is known for many pairs $(G,K)$.
The second ingredient, involves $\Gamma$, and in practice, to determine $\dim V_\pi^\Gamma$ for all $\pi\in\widehat G_\tau$ can be very difficult, even for very simple choices of $G$, $K$ and $\Gamma$.

\subsection{Quotients by cyclic subgroups}
We first recall some basic definitions on the representation theory of compact Lie
groups (see \cite{Knapp-book-beyond} for details).

Let $T$ be a maximal torus of $G$ with Lie algebra $\mathfrak t$.
The complexification $\mathfrak t_\C$ of $\mathfrak t$ is a Cartan subalgebra of $\mathfrak g_\C$.
The associated root system will be denoted by $\Phi(\mathfrak g_\C,\mathfrak t_\C)$.
If $\alpha\in\Phi(\mathfrak g_\C,\mathfrak t_\C)$, then $(\mathfrak g_\C)_\alpha:= \{Y\in\mathfrak g_\C: [X,Y]= \alpha(X) Y \text{ for all }X\in\mathfrak t_\C\}\neq0$.

Given $\pi:G\to \GL(V_\pi)$ a complex finite dimensional representation of $G$, an element $\mu\in \mathfrak t_\C^*$ is called a \emph{weight of $\pi$} if
\begin{equation}
V_\pi(\mu):= \{v\in V_\pi : \pi(\exp(X)) \cdot v = e^{\mu(X)}\, v \text{ for all }X\in\mathfrak t\}\neq0.
\end{equation}
The \emph{weight lattice $P(G)$ of $G$} is the lattice generated by the weights of all finite dimensional representations of $G$.
The weight decomposition of $\pi$ is given by
\begin{equation}
V_\pi = \bigoplus_{\mu \in P(G)} V_\pi(\mu),
\end{equation}
with $V_\pi(\mu)=0$ for all but finitely many $\mu\in P(G)$.
The number $m_{\pi}(\mu):=\dim V_\pi(\mu)$ is the \emph{multiplicity of the weight $\mu$ in $\pi$}.

In what follows we will always assume that $\Gamma$ is cyclic.
By conjugating $\Gamma$ by $x\in G$, if necessary,
we may assume that $\Gamma$ is a subgroup of $T$ without changing the spectrum, since  in formula \eqref{eq:multeigenvalue}, $\dim V_\pi^\Gamma=\dim V_\pi^{x^{-1}\Gamma x}$ for every $x\in G$.

In fact, we may relax the condition to $\Gamma$ being a finite subgroup of $T$ (which is necessarily abelian) without changing anything in the sequel.

Since $\Gamma\subset T$, every element of $\Gamma$ acts diagonally on $V_\pi$ with respect to the weight decomposition.
More precisely, if $v=\sum_{\mu\in P(G)} v_\mu $ with $v_\mu\in V_\pi(\mu)$ for all $\mu\in P(G)$, we have that $\pi(\gamma)\cdot v = \sum_{\mu\in P(G)}  \gamma^\mu v_\mu$, where  $\gamma^\mu=e^{\mu(X_\gamma)}$ for any element $X_\gamma$ in $\mathfrak t$ satisfying that $\gamma=\exp(X_\gamma)$.
It follows that $v\in V_\pi^\Gamma$ if and only if $v_\mu\in V_\pi(\mu)^\Gamma$ for all $\mu\in P(G)$.
We deduce that
\begin{equation}\label{eq:dimV_pi^Gamma}
\dim V_\pi^\Gamma = \sum_{\mu \in P(G)} \dim V_\pi(\mu)^\Gamma
=\sum_{\mu \in \mathcal L_\Gamma} \dim V_\pi(\mu)
=\sum_{\mu \in \mathcal L_\Gamma} m_{\pi}(\mu) ,
\end{equation}
where
\begin{equation}\label{eq:mathcalL_Gamma}
\mathcal L_\Gamma=\{\mu\in P(G): \gamma^\mu=1\;\text{ for all }\gamma\in \Gamma \}.
\end{equation}
It is easy to see, since $\Gamma\subset T$ is a finite group, that $\mathcal L_\Gamma$ is a full lattice in $P(G)$.

We conclude that, when $\Gamma\subset T$, the second ingredient to determine the spectrum of $\Delta_{\tau,\Gamma}$ has been reduced to obtain
an explicit expression for the weight multiplicities $m_{\pi}(\mu)$ of all $\pi\in\widehat G_\tau$ and $\mu\in \mathcal L_\Gamma$.

\subsection{Spectra on functions of lens spaces} \label{subsec:alternative}
The goal of this subsection is to recover the results from Section~\ref{sec:spec0lentes} by using the general approach introduced in the previous subsections.

We set $G=\SO(2n)$ and $K=\{g\in G: g\cdot e_{2n}=e_{2n}\}\simeq \SO(2n-1)$,
thus $G/K$ is diffeomorphic to the $(2n-1)$-dimensional sphere.
We pick the $\Ad(G)$-invariant inner product on $\mathfrak g$ given by $\langle X,Y\rangle =  -\frac12 \operatorname{Tr}(XY)$, which gives constant sectional curvature one.
In other words, the corresponding Riemannian homogeneous space $(G/K, \langle\cdot,\cdot\rangle)$ is isometric to the unit sphere $S^{2n-1}$ in $\R^{2n}$.

We pick $\tau$ the trivial representation $\tau_0$ of $K$, so $E_{\tau_0}$ is the trivial $1$-dimensional vector bundle $S^{2n-1}\times \C$ whose sections are just (complex-valued) functions.
It is well known that $\widehat G_{\tau_0}$ is the collection $\{(\pi_k,\HH_k)\}_{k\geq0}$ (the \emph{spherical representations}), where $\HH_k$ was defined in \eqref{eq:PP_kHH_k}  (but now with $n$ replaced by $2n$) and the action of $g\in G$ on $f\in\HH_k$ is given by
\begin{equation}
(\pi_k(g)\cdot f) (x) = f(g^{-1}x)\qquad (x=(x_1,\dots,x_{2n})^t).
\end{equation}
Furthermore, $\dim \Hom_K(\HH_k,\C)=1$ for all $k\geq0$, and $\lambda(C,\pi_k)=k(k+2n-2)$.

Let $\Gamma$ be any finite subgroup of $G$.
So far, \eqref{eq:multeigenvalue} yields that any eigenvalue of the Laplace--Beltrami operator on $\Gamma\ba S^{2n-1}$ is of the form $k(k+2n-2)$ for some $k\geq0$, and the multiplicity of $k(k+2n-2)$ equals $\dim \HH_k^\Gamma$, which is exactly \eqref{eq:specGammabaG}.

We pick the maximal torus of $G=\SO(2n)$ given by
\begin{equation}
T:=\left\{t(\theta_1,\dots,\theta_n):=
\diag\left(
\left[\begin{smallmatrix}\cos(2\pi\theta_1)& \sin(2\pi\theta_1) \\ -\sin(2\pi\theta_1)&\cos(2\pi\theta_1)
\end{smallmatrix}\right]
,\dots,
\left[\begin{smallmatrix}\cos(2\pi\theta_n)& \sin(2\pi\theta_n) \\ -\sin(2\pi\theta_n)&\cos(2\pi\theta_n)
\end{smallmatrix}\right]
\right)
:\theta_1,\dots,\theta_n \in\R
\right\}.
\end{equation}
The Lie algebra of $T$ is given by
\begin{equation}\label{eq2:h_0}
\mathfrak t =\left\{ X(\theta_1,\dots,\theta_n) :=
\diag\left(
\left[\begin{smallmatrix}0&-2\pi \theta_1\\ 2\pi \theta_1&0\end{smallmatrix}\right]
, \dots,
\left[\begin{smallmatrix}0&-2\pi \theta_n\\ 2\pi \theta_n&0\end{smallmatrix}\right]
\right)
:\theta_1,\dots,\theta_n \in\R
\right\}.
\end{equation}
We identify the complexification $\mathfrak t_\C$ with the set of elements $X(\theta_1,\dots,\theta_n)$ with $\theta_1,\dots,\theta_n\in\C$.
We let $\varepsilon_j\in\mathfrak t_\C^*$, given by $\varepsilon_j(X(\theta_1,\dots,\theta_n))=2\pi \mi\theta_j$ for any $1\leq j\leq n$.
It turns out that the roots are of the form $\pm\varepsilon_i\pm\varepsilon_j$ for $1\leq i<j\leq n$.
Furthermore, the weight lattice $P(G)$ is given by integer combinations of $\varepsilon_1,\dots,\varepsilon_n$, that is, $P(G)=\bigoplus_{j=1}^n\Z\varepsilon_j$.

We now assume $\Gamma=\langle \gamma\rangle $, with $\gamma$ as in \eqref{eq:gamma}, thus $\Gamma\ba S^{2n-1}$ is the lens space $L(q;s_1,\dots,s_n)$.
Note that $\gamma=t(\frac{s_1}q,\dots,\frac{s_n}q)$.
Let $\mu=\sum_j a_j\varepsilon_j$, with $a_j\in\Z$.
Clearly, $\gamma^\mu={\exp(X(\frac{s_1}q,\dots,\frac{s_n}q))}^\mu = e^{\frac{2\pi i}q\sum_j a_j {s_j}} = 1$ if and only if $\sum_j a_j s_j\equiv 0 \pmod q$.
Hence, in this case, as expected, $\mathcal L_\Gamma=\mathcal L(q;s_1,\dots,s_n)$ (see \eqref{eq:congruence-lattice}),  under the identification $P(G)\equiv \Z^n$ given by $\sum_j a_j\varepsilon_j \leftrightarrow (a_1,\dots,a_n)$.

Now, \eqref{eq:dimV_pi^Gamma} yields that the multiplicity of the eigenvalue $\lambda_k=k(k+2n-2)$ is equal to
\begin{equation}\label{eq:dimHH_k^Gamma2}
\dim \HH_k^\Gamma
= \sum_{\mu\in\mathcal L_\Gamma} m_{\pi_k}(\mu)
= \sum_{(a_1,\dots,a_n)\in\Z^n:\;  q\mid a_1s_1+\dots+a_ns_n}  m_{\pi_k}(a_1\varepsilon_1 +\dots+  a_n\varepsilon_n).
\end{equation}
Therefore, an explicit expression for the weight multiplicities of all spherical representations  is sufficient to explicitly describe the spectrum of the lens space $L(q;s_1,\dots,s_n)$.
In fact, by a straightforward computation (see \cite[Lem.~3.2]{LMR-onenorm}), one obtains that
\begin{equation}\label{eq:multiplicity}
m_{\pi_k}(\mu) =
\begin{cases}
\binom{r+n-2}{n-2} & \quad\text{if } \norma\mu=k-2r \text{ for some }r\in \Z_{\geq0},\\
0& \quad\text{otherwise,}
\end{cases}
\end{equation}
where $\norma{\mu}=|a_1|+\dots+|a_n|$ if $\mu=a_1\varepsilon_1 +\dots+  a_n\varepsilon_n\in P(G)$.
It is important to note that $m_{\pi_k}(\mu)$ depends only on $\norma{\mu}$.

Substituting \eqref{eq:multiplicity} into \eqref{eq:dimHH_k^Gamma2} gives
\begin{equation}\label{eq:dimHH_k^Gamma3}
\dim \HH_k^\Gamma
= \sum_{r=0}^{\lfloor k/2\rfloor } \sum_{\mu\in \mathcal L_{\Gamma}:\; \norma{\mu}=k-2r} m_{\pi_k} (\mu)
= \sum_{r=0}^{\lfloor k/2\rfloor }  \tbinom{r+n-2}{n-2} N_{\mathcal L_\Gamma}(k-2r),
\end{equation}
which is the same conclusion as in \eqref{eq:dimHH_k^Gamma}.
By \eqref{eq:F_Gamma(z)pto} in the reverse direction, it follows that
\begin{equation}\label{eq:F_L(z)2}
 F_{\Gamma}(z) = \frac{\vartheta_{\mathcal L_{\Gamma}}(z) }{(1-z^2)^{n-1}},
\end{equation}
which is precisely \eqref{eq:F_L(z)}.

\section{Hodge spectra of lens spaces} \label{sec:p-spectrum}

In this section we study the spectrum of the Hodge--Laplace operator acting on $p$-forms of a lens space.

\subsection{Hodge spectra of spherical space forms} \label{subsec:p-spec-arbitrario}

Let $M$ be a compact Riemannian manifold.
A $p$-form is a smooth section of the $p$-exterior power of the cotangent bundle (i.e.\ $\bigwedge^p T^*M$).
The Hodge--Laplace operator acts on smooth complex $p$-forms as $dd^*+d^*d$, extending the Laplace--Beltrami operator on $0$-forms (i.e.\ functions).
The spectrum of this operator is called the \emph{$p$-spectrum of $M$}.
Furthermore, two Riemannian manifolds or orbifolds with the same $p$-spectra are said to be \emph{$p$-isospectral}.

Throughout this section, we let $G=\SO(2n)$, $K=\{g\in G: g\cdot e_{2n}=e_{2n}\}\simeq \SO(2n-1)$, $G/K\simeq S^{2n-1}$ and $\Gamma$ a finite subgroup of $G$.

For $0\leq p\leq 2n-1$, we consider the operator $\Delta_p$ given by the Hodge--Laplace operator acting on $\Gamma$-invariant $p$-forms on $S^{2n-1}$. 
When $\Gamma$ acts freely on $S^{2n-1}$, $\Delta_p$ is precisely the Hodge--Laplace operator on $\Gamma\ba S^{2n-1}$  the $p$-forms on $\Gamma\ba S^{2n-1}$ are identified with the $\Gamma$-invariant $p$-forms on $S^{2n-1}$. 
In light of the discussion on manifolds and orbifolds in Section~\ref{sec:spectraspherical}, we shall call  $\Delta_p$ 
the Hodge--Laplace operator on $p$-forms on $\Gamma\ba S^{2n-1}$,  even in the case when $\Gamma\ba S^{2n-1}$ has singularities.
We will assume $p\leq n-1$, given that the spectra on $p$-forms and on $(2n-1-p)$-forms coincide, since $\Gamma\ba S^{2n-1}$ is orientable.

The bundle $\bigwedge^p T^*(\Gamma\ba S^{2n-1})$  is isomorphic to  $\Gamma\ba E_{\tau_p}$ (as defined in  Subsection~\ref{subsec:E_tau}), where $\tau_p$ is the   $p$-exterior power of the complexified  standard representation of $K$ on $\bigwedge^p (\C^{2n-1})$
and the operator $\Delta_{\tau_p,\Gamma}$ is identified with the Hodge--Laplace operator acting on $p$-forms on $\Gamma\ba S^{2n-1}$  (see proof in  \cite[Prop.~2.3]{IkedaTaniguchi78} for any compact symmetric space).
The next step will be to describe the representations in $\widehat G_{\tau_p}$.

We keep the root system $\Phi(\mathfrak g_\C,\mathfrak t_\C)$ from Subsection~\ref{subsec:alternative}, and the lexicographic order on $\mathfrak t_\C^*$ with respect to the basis $\{\varepsilon_1,\dots,\varepsilon_n\}$.
The corresponding positive roots are $\Phi^+(\mathfrak g_\C,\mathfrak t_\C):=\{\varepsilon_i\pm \varepsilon_j: i<j\}$, $\rho=\sum_{j=1}^n (n-j)\varepsilon_j$ and
furthermore, the set of dominant weights in $P(G)$ is
\begin{equation}
P^{++}(G):= \{a_1\varepsilon_1+\dots+a_n\varepsilon_n \in P(G): a_1\geq a_2\geq\dots a_{n-1}\geq |a_n|\}.
\end{equation}

The irreducible representations of $G$ are parametrized by the elements in $P^{++}(G)$, by the highest weight theorem.
For instance, the highest weight of $(\pi_k,\HH_k)$ is $k\varepsilon_1$.
Note that \eqref{eq:Casimireigenvalue} gives the eigenvalue
\begin{equation}
\lambda_k=\langle k\varepsilon_1,k\varepsilon_1+2\rho\rangle = k(k+2n-2).
\end{equation}

The classical branching rule from $G$ to $K$ (see for instance \cite[Thm.~9.16]{Knapp-book-beyond}) gives an explicit parametrization of the irreducible representations of $G$ occurring in $\widehat G_{\tau_p}$ for every $p$.
Namely,
\begin{equation}
\widehat G_{\tau_p} =
\begin{cases}
\{1_G\}\cup \{\pi_{k,1}:k\geq 0\}
	\quad&\text{ if }p=0,\\
\{\pi_{k,p},\pi_{k,p+1}: k\geq0\}
	\quad&\text{ if }1\leq p\leq n-2,\\
\{\pi_{k,n-1}, \pi_{k,n}^+,\pi_{k,n}^-: k\geq0\}
	\quad&\text{ if }p=n-1,
\end{cases}
\end{equation}
where $\pi_{k,p}$ has highest weight $\Lambda_{k,p}:= k\varepsilon_1  + \sum_{j=1}^p \varepsilon_j$ for every $k\geq0$ and $1\leq p \leq n-1$, and the representation $\pi_{k,n}^\pm$ has highest weight $\Lambda_{k,n}^\pm:= k\varepsilon_1  + \sum_{j=1}^{n-1} \varepsilon_j \pm \varepsilon_n$ for every $k\geq0$.
 Moreover, the branching law from $G$ to $K$ is \emph{multiplicity free}, i.e.\ $\dim\Hom_K(\tau_p,\pi)\leq1$ for every $\pi\in\widehat G$, thus $\tau_p$ occurs exactly once in the decomposition of $\pi|_K$ in irreducible constituents for every $\pi\in\widehat G_{\tau_p}$.

The eigenvalue corresponding to the representation $\pi_{k,p}$ (or $\pi_{k,n}^\pm$ in case $p=n$) is given by
\begin{align}\label{eq:lambda_kp}
\lambda_{k,p}&:= \langle \Lambda_{k,p},\Lambda_{k,p}+2\rho \rangle
= (k+1)(k+1+2n-2) + \textstyle\sum\limits_{j=2}^p (1+2n-2j) \\
&= k^2+2nk+(2n-1)+ (2n+1)(p-1)-p(p+1)+2\notag\\
&= (k+p)(k+2n-p).  \notag
\end{align}
We set $\pi_{k,0}=1_G$, $\lambda_{k,0}=0$, and $\pi_{k,n}=\pi_{k,n}^+\oplus \pi_{k,n}^-$ for any $k\geq0$.
It is easy to see that $\lambda_{k,p}\neq \lambda_{k',p+1}$ for all $k,k'\geq0$ and $0\leq p\leq n-1$.
 This new notation is related to the notation in Subsection~\ref{subsec:alternative} as follows: $\pi_{k,1}=\pi_{k+1}$ and $\lambda_{k,1}=\lambda_{k+1}$.

Collecting this information, the $p$-spectrum of $S^{2n-1}$ is given in \cite{IkedaTaniguchi78}.
Thus, \eqref{eq:multeigenvalue} describes the $p$-spectrum of $\Gamma\ba S^{2n-1}$ as follows.

\begin{proposition}\label{prop:p-spec}
Let $\Gamma$ be a finite subgroup of $G$ and $0\leq p\leq n-1$.
Any eigenvalue in the spectrum of $\Delta_{\tau_p,\Gamma}$ is of the form $\lambda_{k,p}$ or $\lambda_{k,p+1}$ for some $k\geq0$, with multiplicity given by $\dim V_{\pi_{k,p}}^\Gamma$ and $\dim V_{\pi_{k,p+1}}^\Gamma$ respectively.
\end{proposition}

\subsection{The spectra of lens spaces}\label{subsec:characterization}
Our next main goal is to obtain a characterization of pairs of lens spaces that are either $0$-isospectral or  $p$-isospectral for every  $p$, in terms of geometric properties of their associated lattices.

We now assume that $\Gamma=\langle \gamma\rangle$ with $\gamma$ as in \eqref{eq:gamma}, thus $\Gamma\ba S^{2n-1}=L(q;s_1,\dots,s_n)$ is a lens space.
Again, we may relax the condition on $\Gamma$ to be a finite subgroup of the maximal torus $T$ of $G$, without many further changes in the sequel.

We abbreviate $L=L(q;s_1,\dots,s_n)$ and $\mathcal L= \mathcal L(q;s_1,\dots,s_n)$.
We recall that $\mathcal L$ is given by elements $\mu= \sum_{j=1}^n a_j\varepsilon_j \in P(G)$ satisfying that $a_1,\dots,a_n\in\Z$ and $a_1s_1+\dots+a_ns_n\equiv 0 \pmod q$.
Similarly as in \eqref{eq:dimHH_k^Gamma2}, \eqref{eq:dimV_pi^Gamma} gives
\begin{equation}\label{eq:dimV_kp^Gamma}
\dim V_{\pi_{k,p}}^\Gamma
= \sum_{\mu\in\mathcal L_\Gamma} m_{\pi_{k,p}}(\mu)
= \sum_{(a_1,\dots,a_n)\in\Z^n:\;  q\mid a_1s_1+\dots+a_ns_n}  m_{\pi_{k,p}}(a_1\varepsilon_1 +\dots+  a_n\varepsilon_n).
\end{equation}

Thus we can see ---as it was shown in detail in \cite[Lem.~3.3]{LMR-onenorm}--- that two weights with the same one-norm and the same number of zero coordinates have the same multiplicity in any $\pi_{k,p}$.
More precisely, if
$Z(\mu):=\#\{i: 1\leq i\leq n,\, a_i=0\}$ for $\mu=\sum_{i=1}^n a_i\varepsilon_i$, we have
\begin{lemma}\label{lem:3.3}  If $\norma{\mu}=\norma{\mu'}$ and  $Z(\mu)=Z(\mu')$
then
\begin{equation*}
m_{\pi_{k,p}}(\mu)=m_{\pi_{k,p}}(\mu'), \quad\text{ for all $k\geq0$ and $1\leq p\leq n$.}
\end{equation*}
\end{lemma}

We now set
\begin{equation}
N_{\mathcal L}(k,\ell) = \#\{\mu\in\mathcal L: \norma{\mu}=k,\, Z(\mu)= \ell\}.
\end{equation}
 By the previous lemma, \eqref{eq:dimV_kp^Gamma}, and the fact that $m_{\pi_{k,p}}(\mu)=0$ for all $\mu\in P(G)$ satisfying $\norma{\mu}>\norma{\Lambda_{k,p}}=k+p$, it follows that
\begin{align}\label{eq:dimV_kp^Gamma2}
\dim V_{\pi_{k,p}}^{\Gamma}
&=\sum_{r=0}^{[(k+p)/2]} \;\sum_{\zz =0}^n \; \sum_{\mu\in\mathcal L: \, \norma{\mu}=k+p-2r,\, Z(\mu)=\ell} m_{\pi_{k,p}}(\mu),
\\
&=\sum_{r=0}^{[(k+p)/2]} \;\sum_{\zz =0}^n \; m_{\pi_{k,p}}(\mu_{r,\zz })\; N_{\mathcal L}(k+p-2r,\zz ),
\notag
\end{align}
where $\mu_{r,\zz }$ is any weight in $\mathcal L$ with $\norma{\mu_{r,\zz }} = k+p-2r$ and having $\zz$ zero coordinates.

Given two lattices, we say that they are $\norma{\cdot}^*$-isospectral if
for each $k\in\N$ and $0\leq \zz\leq n$ there are the same number of elements $\mu$ in each lattice having $\norma{\mu} = k$ and having exactly $\zz$ coordinates equal to zero, i.e.\ $N_{\mathcal L}(k,\ell)=N_{\mathcal L}(k,\ell)$ for all $k\geq0$ and $0\leq \ell\leq n$.

As an output of \eqref{eq:dimV_kp^Gamma2}, we prove one of the main results in \cite{LMR-onenorm}
(Thm.~3.6).
We note that for the converse assertion, we need to invert \eqref{eq:dimV_kp^Gamma2},
showing that one can solve for the $N_{\mathcal L}(k+p-2r,\zz)$ in terms of the values $\dim V_{\pi_{k,p}}^{\Gamma}$.
We include in the statement part (i) for completeness.

\begin{theorem} \label{thm:3.6(ii)}
Let $L=\Gamma\ba S^{2n-1}$ and $L'=\Gamma'\ba S^{2n-1}$ be lens spaces with associated congruence lattices $\mathcal L$ and $\mathcal L'$ respectively.
Then
	\begin{enumerate}
		\item[(i)] $L$ and $L'$ are $0$-isospectral if and only if $\mathcal L$ and $\mathcal L'$ are $\norma{\cdot}$-isospectral.
		\item[(ii)] $L$ and $L'$ are $p$-isospectral for all $p$ if and only if $\mathcal L$ and $\mathcal L'$ are    $\norma{\cdot}^*$-isospectral.
	\end{enumerate}
\end{theorem}

\begin{remark}\label{rem:isometria}
We recall the classical well-known fact (see \cite[Ch.~V]{Cohen:book}) that two lens spaces $L=L(q;s_1,\dots,s_n)$ and $L'=L(q;s_1',\dots,s_n')$ are isometric if and only if
\begin{itemize}
\item[($*$)] There exist $\sigma$ a permutation of $\{1,\dots,n\}$, $\epsilon_1,\dots,\epsilon_n\in\{\pm1\}$ and $t\in\Z$ coprime to $q$ such that $s_{\sigma(j)}'\equiv t\epsilon_js_j\pmod{q}$
\textrm{ for every }  $1\leq j\leq n.$
\end{itemize}

We observe that we may restate condition ($*$)  in terms of the geometry of the associated congruence lattices as follows (see \cite[Prop. 3.3]{LMR-onenorm}):
\begin{enumerate}
\item[($*'$)] There is a linear one-norm isometry from $\mathcal L(q;s_1,\dots,s_n)$ onto $\mathcal L(q;s_1',\dots,s_n')$.
\end{enumerate}

Similarly, we may rephrase the assertion in Theorem~\ref{thm:3.6(ii)} (i) by saying that two lens spaces are isospectral if and only if the associated lattices are one-norm isometric (but not necessarily via a linear map).
\end{remark}

\begin{remark}\label{rem3:ej-Ikeda}
Ikeda in \cite{Ikeda80_isosp-lens} gave many pairs of non-isometric lens spaces that are $0$-isospectral.
The simplest such pair is $L(11;1,2,3)$ and $L(11;1,2,4)$ in dimension $5$.
In light of Theorem~\ref{thm:3.6(i)}, $\mathcal L=\mathcal L(11;1,2,3)$ and $\mathcal L'=\mathcal L(11;1,2,4)$ must be $\norma{\cdot}$-isospectral (i.e.\ $N_{\mathcal L}(k)= N_{\mathcal L'}(k)$ for all $k\geq0$).
However, it is a simple matter to check that $\mathcal L$ and $\mathcal L'$ are not $\norma{\cdot}^*$-isospectral since $N_{\mathcal L}(3,0)=2 \neq N_{\mathcal L'}(3,0)=0$ and  $N_{\mathcal L}(3,1)=2 \neq N_{\mathcal L'}(3,1)=4$, contradicting the condition in Theorem~\ref{thm:3.6(ii)}.
	
Ikeda and Yamamoto showed that two $0$-isospectral $3$-dimensional lens spaces are isometric (\cite{IkedaYamamoto79}, \cite{Yamamoto80}).
Also, in the relevant paper \cite{Ikeda88}, Ikeda produced for each given $p_0$ pairs of lens spaces that are $p$-isospectral for every $0\le p \le p_0$ but are not $p_0+1$ isospectral.
We refer the reader to \cite[\S4]{LMR-survey} for a detailed account on these results.
\end{remark}

\subsection{Non strongly isospectral manifolds $p$-isospectral for all $p$} \label{subsec:examples}

Two compact Riemannian manifolds are said to be \emph{strongly isospectral} if they are isospectral with respect to every natural strongly elliptic operator acting on sections of a natural vector bundle (see \cite{DeT-Gordon}).
In particular, strongly isospectral manifolds are $p$-isospectral for all $p$.
The famous Sunada method always produces strongly isospectral manifolds.
Most of the isospectral pairs known are also strongly isospectral,
however, along the years, a number of isospectral examples have been constructed that do not follow the Sunada construction.
For instance, there are examples isospectral on functions but not on $1$-forms (see \cite{Gordon86}),
lens spaces $p$-isospectral for  $0\le p\le p_0$ but not $(p_0+1)$-isospectral (see \cite{Ikeda88}) and also many in the case of flat manifolds (\cite{MRpiso}).

The first examples of non-isometric compact connected Riemannian manifolds $p$-isospectral for all $p$, but still not strongly isospectral were constructed in \cite{LMR-onenorm}, where we exhibit infinitely many pairs of $\norma{\cdot}^*$-isospectral congruence lattices.

\begin{theorem}
For any $r\geq7$ and $t$ positive integers such that $r$ is coprime to $3$, the 3-dimensional congruence lattices
\begin{equation}\label{eq:LMRexamples}
\mathcal L(r^2t;\;1,\;1+rt,\;1+3rt)
\quad\text{and}\quad
\mathcal L(r^2t;\;1,\;1-rt,\;1-3rt)
\end{equation}
are $\norma{\cdot}^*$-isospectral.
\end{theorem}

This result implies the corresponding result for lens spaces, thus we produce an infinite family of pairs of $5$-dimensional lens spaces that are $p$-isospectral for all $p$. As a consequence, by using a result of Ikeda, we obtain pairs of lens spaces that are $p$-isospectral for all $p$ in arbitrarily high dimensions, and also pairs of compact Riemannian manifolds with this property in every dimension $n\ge 5$.
We point out that the resulting lens spaces are homotopically equivalent but not homeomorphic to each other (see \cite[Lemma~7.1]{LMR-onenorm}).

To prove that the corresponding lens spaces for $r\geq 7$ are not strongly isospectral, we first use that for arbitrary $\Gamma,\Gamma'$ finite subgroups of $\SO(2n)$, the spaces $\Gamma\ba S^{2n-1}$ and $\Gamma'\ba S^{2n-1}$ are strongly isospectral if and only if $\Gamma$ and $\Gamma'$ are almost conjugate in $\Ot(2n)$.
 This follows from Pesce's characterization of strongly isospectral spherical space forms (\cite[Prop.~III.1]{Pesce95}). 
  The assertion that $\Gamma$ and $\Gamma'$ are representation equivalent in $\Ot(2n)$ is equivalent to be almost conjugate in $\Ot(2n)$ (see for instance \cite[Lem.~ 2.12]{Wolf01}).
Since lens spaces are quotients of $S^{2n-1}$ by cyclic subgroups, and almost conjugate cyclic subgroups are necessarily conjugate, then two lens spaces are strongly isospectral if and only if they are isometric. 
In \cite[Prop~5.4]{LMR-onenorm} it is shown that the  congruence lattices in \eqref{eq:LMRexamples} are not $\norma\cdot$-isometric, and therefore the corresponding lens spaces are not isometric by Remark~\ref{rem:isometria}.

\begin{remark}
In \cite{LMR-onenorm}, a finiteness condition for $\norma{\cdot}^*$-isospectrality of two congruence lattices was given, by counting the number of lattice points of a fixed norm $k$ in a small cube, having exactly $\zz$ zero coordinates.
This allowed the authors to give tables with many examples of low-dimensional pairs of
lens spaces $p$-isospectral for all $p$ (see \cite[Tables~1--2]{LMR-onenorm} and \cite[Tables 1--3]{LMR-survey}).
\end{remark}

\begin{remark}
DeFord and Doyle~\cite{DeFordDoyle14} made a big step toward a full classification of lens spaces $p$-isospectral for all $p$.
In particular, for fixed non-negative integers $r$ and $t$, they found a sufficient condition on $(d_1,\dots,d_{n})\in\Z^{n}$ such that
\begin{equation}\label{eq:LMR(r,t,d)}
L(r^2t; 1+d_1rt,\dots, 1+d_{n}rt )
\quad\text{and}\quad
L(r^2t; 1-d_1rt,\dots, 1-d_{n}rt )
\end{equation}
are $p$-isospectral for all $p$.
Indeed, they showed that  it suffices that, for each positive divisor $u$ of $r$, either the numbers $d_1,\dots,d_n$ modulo $u$ are all different, or $L(u^2t; 1+d_1ut,\dots, 1+d_{n}ut )$ and $L(u^2t; 1-d_1ut,\dots, 1-d_{n}ut )$ are isometric.

For instance, it is easy to see that $(d_1,d_2,d_3)=(0,1,3)$ satisfies DeFord--Doyle's condition when $r$ is not divisible by $3$.
This shows that the pair of lens spaces in \eqref{eq:LMRexamples} are $p$-isospectral for all $p$.

Curiously, every example found by the computer so far satisfies this condition, with the only exception of the pair $L(72; 1, 5, 7, 17, 35)$ and $L(72; 1, 5, 7, 19, 35)$.
The pair $L(72; 1, 5, 7, 11, 19, 25, 35)$ and $L(72; 1, 5, 7, 11, 23, 29, 31)$, which is dual to the previous one in a certain sense (see \cite[Thm.~7.3]{LMR-onenorm} and \cite[\S5.3]{LMR-survey}), does not fit in DeFord--Doyle's method either.
\end{remark}

We point out that all the known families of lens spaces $p$-isospectral for all $p$ have only two elements, which
suggests the following question
\begin{problem}\cite[Ques.~5.3]{LMR-onenorm}
\cite[Ques.~5.8]{LMR-survey}
Are there families of non-isometric lens spaces $p$-isospectral for all $p$ having at least three elements?
\end{problem}

The next problem is a major one in the spectral theory of lens spaces (see also \cite[\S12(5)]{DeFordDoyle14}).
\begin{problem}\cite[Prob.~5.7]{LMR-survey}
Determine all families of $n$-dimensional lens spaces $p$-isospectral for all $p$, with fundamental group of order $q$, for every $q$ and $n$.
\end{problem}

\subsection{Spectra on $p$-forms of lens spaces for an individual $p$.}

If we wish to determine conditions for the $p$-isospectrality of two lens space $L$, $L'$, for \emph{one} individual $p\ne 0$, formula \eqref{eq:dimV_kp^Gamma2} does not allow to obtain a neat geometric condition on the corresponding congruence lattices $\mathcal L$, $\mathcal L'$, as was the case for $p=0$ or for all $p$'s simultaneously.

Nevertheless, Rossi Bertone and the first named author were able to give a closed explicit formula for the multiplicity $m_{\pi_{k,p}}(\mu)$ in terms of $\norma{\mu}$ and $Z(\mu)$ (\cite[Thm.~IV.1]{LR-fundstring}).
Such expression is pretty long, unlike the neat formula \eqref{eq:multiplicity},
however, it was sufficient to give in \cite{Lauret-p-spectralens} a refined description of every individual $p$-spectrum of a lens space in place of the expression in \eqref{eq:dimV_kp^Gamma2}.
To state the result, we first introduce some notation.

Extending \eqref{eq:spectralgeneratingfunction}, Ikeda in \cite{Ikeda88} associated to a quotient $\Gamma\ba S^{2n-1}$ ($\Gamma$ an arbitrary finite subgroup of $G=\SO(2n)$) the functions
\begin{equation}\label{eq:F_Gamma^p}
F_\Gamma^{p}(z) =\sum_{k\geq0} \dim V_{\pi_{k,p+1}}^\Gamma z^k,
\qquad\text{for $0\leq p\leq n-1$.}
\end{equation}
(Note that $F_{\Gamma}(z) = 1+zF_{\Gamma}^0(z)$.)
Its $k$-th term is the multiplicity of the eigenvalue $\lambda_{k,p+1}$ (defined in \eqref{eq:lambda_kp}) in the $p$-spectrum and also in the $(p+1)$-spectrum of $\Gamma\ba S^{2n-1}$ (see Proposition~\ref{prop:p-spec}).
It follows immediately that
\begin{equation}\label{eq:Ikeda'scharacterizationp-iso}
\Gamma\ba S^{2n-1}\text{ and } \Gamma'\ba S^{2n-1}
\text{ are $p$-isospectral }
\Longleftrightarrow \
F_{\Gamma}^{p-1}(z)= F_{\Gamma'}^{p-1}(z) \text{ and } F_{\Gamma}^p(z)= F_{\Gamma'}^p(z).
\end{equation}
Here, it is understood that $F_{\Gamma}^{p-1}(z)=0$ for $p=0$.

The main result in \cite{Lauret-p-spectralens} states, for $\Gamma$ any finite subset of $T$, that
\begin{equation}\label{eq1:F_Gamma^p-1-formula}
F_\Gamma^{p-1}(z) = \frac{1}{(1-z^2)^{n-1}}\sum_{\ell=0}^{n} \vartheta_{\mathcal L_\Gamma}^{(\ell)}(z) \; A_{p}^{(\ell)}(z)+\frac{(-1)^p}{z^p},
\end{equation}
where
\begin{align}
\vartheta_{\mathcal L_\Gamma}^{(\ell)}(z)&= \sum_{k\geq0} N_{\mathcal L_\Gamma}(k,\ell)\, z^k  \qquad\text{(compare with \eqref{eq:one-normtheta}), }\\	
	\label{eq3:A_pl}
	A_{p}^{(\ell)}(z) &= \sum_{j=1}^{p} (-1)^{j-1}  \sum_{t=0}^{\lfloor\frac{p-j}{2}\rfloor} \binom{n-p+j+2t}{t}
	\\ &\quad
	\sum_{\beta=0}^{p-j-2t} 2^{p-j-2t-\beta} \binom{n-\ell}{\beta} \binom{\ell}{p-j-2t-\beta}
	\sum_{\alpha=0}^\beta \binom{\beta}{\alpha} \sum_{i=0}^{j-1} z^{p-2(j+t+\alpha-i)}. \notag
\end{align}
In conclusion, through long calculations, one can explicitly compute the multiplicity of any eigenvalue in the $p$-spectrum of $\Gamma\ba S^{2n-1}$ in terms of the numbers $\{N_{\mathcal L_\Gamma}(k,\ell): k\geq0,\, 0\leq \ell\leq n\}$.

Furthermore, $\vartheta_{\mathcal L_\Gamma}^{(\ell)}(z)$ is a rational function with a common denominator $(1-z^q)^{n-\ell}$,  where $q$ stands for the number of elements of $\Gamma$ (see \cite[Thm.~2.4]{Lauret-p-spectralens}).
This implies that $F_{\Gamma}^{p-1}(z)$ can be determined with finitely many computations, and consequently, one has a finiteness condition for two spaces $\Gamma\ba S^{2n-1}$ and $\Gamma'\ba S^{2n-1}$ to be $p$-isospectral for a single $p$.
This allowed the computational study \cite{Lauret-computationalstudy} of $p$-isospectral lens spaces, which will be summarized in Subsection~\ref{subsec:computationalstudy}.

Ikeda's characterization \eqref{eq:Ikeda'scharacterizationp-iso} of $p$-isospectral spherical space forms immediately implies, for a fixed $0\leq p_0\leq n-1$, that
$\Gamma\ba S^{2n-1}$ and $\Gamma'\ba S^{2n-1}$ are $p$-isospectral for all $0\leq p\leq p_0$ if and only if $
F_{\Gamma}^{p}(z)= F_{\Gamma'}^{p}(z)$ for all $0\leq p\leq p_0$.
This condition was used by Ikeda in \cite{Ikeda88} to construct, for any given $p_0$, examples of lens spaces $p$-isospectral for all $0\leq p\leq p_0$ but not $(p_0+1)$-isospectral (see \cite[\S4]{LMR-survey} for a discussion on this result).

A careful manipulation of the functions $F_{\Gamma}^{0}(z), F_{\Gamma}^{1}(z),\dots, F_{\Gamma}^{p}(z)$ by using the expression \eqref{eq1:F_Gamma^p-1-formula} gave the following geometric characterization of lens spaces $p$-isospectral for all $0\leq p\leq p_0$.

\begin{theorem}\cite[Cor.~2.3]{Lauret-p-spectralens} \label{cor:charact[0,p]}
Let $0\leq p_0\leq n-1$ and let $\Gamma$ and $\Gamma'$ be finite subgroups of the maximal torus $T$ of $G$.
Then, $\Gamma\ba S^{2n-1}$ and $\Gamma'\ba S^{2n-1}$ are $p$-isospectral for all $0\leq p\leq p_0$ if and only if
	\begin{equation}\label{eq1:condition[0,p_0]-isosp}
	\sum_{\ell=0}^n \ell^h \, \vartheta_{\mathcal L_{\Gamma}}^{(\ell)}(z) =
	\sum_{\ell=0}^n \ell^h \, \vartheta_{\mathcal L_{\Gamma'}}^{(\ell)}(z)
	\qquad\text{for all }0\leq h\leq p_0,
	\end{equation}
	or equivalently,
	\begin{equation}\label{eq1:condition[0,p_0]-isospN_kl}
	\sum_{\ell=0}^n \ell^h \, N_{\mathcal L_{\Gamma}}(k,\ell) =
	\sum_{\ell=0}^n \ell^h \, N_{\mathcal L_{\Gamma'}}(k,\ell)
	\quad\text{for all } k\geq0,\qquad\text{for all }0\leq h\leq p_0.
	\end{equation}
\end{theorem}

Clearly, the case $p_0=0$ in this result coincides with Theorem~\ref{thm:3.6(ii)}(i) since $\vartheta_{\mathcal L_{\Gamma}}(z)=\sum_{\ell=0}^n \vartheta_{\mathcal L_{\Gamma}}^{(\ell)}(z)$, and one has that $\vartheta_{\mathcal L_{\Gamma}}(z)=\vartheta_{\mathcal L_{\Gamma'}}(z)$ if and only if $\mathcal L_{\Gamma}$ and $\mathcal L_{\Gamma'}$ are $\norma{\cdot}$-isospectral.
Moreover, it turns out that the condition \eqref{eq1:condition[0,p_0]-isosp} for $p_0=n-1$ is equivalent to $\mathcal L_{\Gamma}$ and $\mathcal L_{\Gamma'}$ being $\norma{\cdot}^*$-isospectral, thus Theorem~\ref{cor:charact[0,p]} for $p_0=n-1$ coincides with Theorem~\ref{thm:3.6(ii)}(ii) (see \cite[Rem.~3.5]{Lauret-p-spectralens}).

\subsection{Some results on $\tau$-spectra for homogeneous vector bundles}
One may consider the homogeneous vector bundle $E_{\tau}$ on $S^{2n-1}$ associated to an arbitrary finite dimensional representation $\tau$ of $K=\SO(2n-1)$.
One usually gets irreducible representations $\pi$ in $\widehat G_\tau$ that are different from $\pi_{k,p}$ for all $k,p$.
The main difficulty  to describe the spectrum of $\Delta_{\tau,\Gamma}$ for any $\Gamma\subset T$ is to compute the multiplicity $m_{\pi}(\mu)$  from \eqref{eq:dimV_pi^Gamma} of an arbitrary weight $\mu$ in $\pi$ in order to determine $\dim V_\pi^\Gamma$.
In general, there are no closed formulas for $m_{\pi}(\mu)$ adequate to be used in \eqref{eq:dimV_pi^Gamma}, neither a visible geometric regularity of the weight multiplicities, like in Lemma~\ref{lem:3.3}.

As a continuation of \cite{LR-fundstring}, in \cite{LR-bivariate} expressions were given for the weight multiplicities of the irreducible representations of $G=\SO(2n)$ having highest weight $k\varepsilon_1+l\varepsilon_2$ for all $k\geq l\geq 0$.
These representations are called \emph{bivariate representations} because their highest weights are integral combinations of the first two fundamental weights.

\begin{problem}
For some $\tau\in \widehat K$ different from every $\tau_p$, use \cite{LR-bivariate} to give an expression for the spectral generating function associated to  $\Delta_{\tau,\Gamma}$, for any $\Gamma\subset T$.
Apply such expression to give a geometric characterization of $\tau$-isospectral lens spaces.
Are there $\tau$-isospectral lens spaces?
\end{problem}

If we choose $K'=\SO(2n-2)\times\SO(2)$, then $G/K'$ is diffeomorphic to the Grassmannian space of $2$-dimensional linear subspaces in $\R^{2n}$.
 The $p$-spectrum of $G/K'$ was studied by  Strese \cite{Str80} for $p=0,1$ and by Tsukamoto~\cite{Tsukamoto81} for any $p$.
If $\tau_0$ denotes the trivial representation of $K'$,
then $\Delta_{\tau_0,\Gamma}$ is the Laplace--Beltrami operator on $\Gamma\ba G/K'$ for any finite subgroup $\Gamma$ of $G$.
Moreover, 
the highest weights of the elements in $\widehat G_{\tau_0}$ are non-negative integral combinations of $2\varepsilon_1$ and $\varepsilon_1+\varepsilon_2$, thus every $\pi$ in $\widehat G_{\tau_0}$ is a bivariate representation. 
Moreover, $\dim \Hom_{K'}(\tau_0,\pi)=1$ for all $\pi\in\widehat G_{\tau_0}$. 
Therefore, according to \eqref{eq:multeigenvalue}, the only remaining ingredient to compute the $0$-spectrum of $\Gamma\ba G/K'$ is to determine $\dim V_\pi^\Gamma$ for every bivariate representation $\pi$.

\begin{problem}\cite[Rem.~III.17]{LR-bivariate} \label{prob:2-Grass}
Describe the spectrum of the Laplace--Beltrami operator on $\Gamma\ba \SO(2n)/(\SO(2n-2)\times\SO(2))$ for any $\Gamma \subset T^n$.
\end{problem}

Ikeda in \cite{Ikeda97} studied the problem whether $0$-isospectrality of the spherical space forms $\Gamma\ba S^{2n-1}$ and $\Gamma'\ba S^{2n-1}$ implies
the $0$-isospectrality of the quotients $\Gamma\ba \SO(2n)/\SO(2n-2m)\times\SO(2m)$ and $\Gamma'\ba \SO(2n)/\SO(2n-2m)\times\SO(2m)$ (of the Grassmannian $\SO(2n)/\SO(2n-2m)\times\SO(2m)$) and
he proved the validity of this assertion for $\Gamma$ and $\Gamma'$ almost conjugate subgroups.
We note that this result follows by an application of Sunada's method: indeed, if $\Gamma$ and $\Gamma'$ are almost conjugate subgroups of a compact Lie group $\widetilde G$, then $\Gamma\ba \widetilde G/\widetilde K$ and $\Gamma'\ba \widetilde G/\widetilde K$ are strongly isospectral for every closed subgroup $\widetilde K$ of $\widetilde G$ (see \cite[Prop.~2.10]{Wolf01}).
For lens spaces $\Gamma\ba S^{2n-1}$ and $\Gamma'\ba S^{2n-1}$ the question remains open.
The case $m=1$ might be tractable assuming some progress in Problem~\ref{prob:2-Grass}.

\begin{problem}
Investigate the existence of $0$-isospectral pairs of the form $\Gamma\ba \SO(2n)/\SO(2n-2)\times\SO(2)$ and $\Gamma'\ba \SO(2n)/\SO(2n-2)\times\SO(2)$ with $\Gamma$ and $\Gamma'$ cyclic subgroups of $\SO(2n)$.
\end{problem}

In a different direction, there is the question whether the lens spaces $L=L(49;1,6,15)$ and $L'=L(49;1,6,20)$ in \eqref{eq:LMRexamples}, that are $p$-isospectral for all $p$, are $\tau$-isospectral for some $\tau\in \widehat K$ different from $\tau_p$.
Several examples of $\tau\in\widehat K$ were given in \cite{LMR-onenorm} such that $L$ and $L'$ are not $\tau$-isospectral.
Actually, this seems to be the case for generic $\tau$.
In \cite[\S8]{LMR-onenorm} the authors claimed the following
\begin{conjecture}
There are only finitely many irreducible representations $\tau$ of $\SO(5)$ such that $L(49;1,6,15)$ and $L(49;1,6,20)$ are $\tau$-isospectral.
\end{conjecture}

\section{Connections with Ehrhart theory and toric varieties} \label{sec:Ehrhart}

In this section, we first give a description of the spectrum of a lens space by using Ehrhart series, as in \cite{MohadesHonari17} (see also \cite[Thm.~3.9]{Lauret-spec0cyclic}), and we end the section by giving a connection with toric varieties, developed by Mohades and Honari in~\cite[\S5]{MohadesHonari17}.

A \emph{(convex) polytope} is the smallest convex set containing a finite set in $\R^n$.
For $\mathcal P$ a polytope and $k$ a positive integer, we denote by $k\mathcal P$ the $k$-th dilate of $\mathcal P$, that is,
$k\mathcal P= \{kx : x\in\mathcal P\}$.
A \emph{rational polytope} (resp.\ \emph{integral polytope}) is a polytope whose vertices have rational (resp.\ integer) coordinates.
Ehrhart~\cite{Ehrhart62} proved that the \emph{integer-point enumerator} $i_{\mathcal P}(k):= \#(k\mathcal P\cap \Z^n)$ of an integral polytope $\mathcal P$
is actually a polynomial in $k$, by studying the associated power series
\begin{equation}
\Ehr_{\mathcal P}(z)
:= \sum_{k\geq0} i_{\mathcal P}(k)\, z^k
=\sum_{k\geq0} \#(k\mathcal P\cap \Z^n)\, z^k,
\end{equation}
known nowadays as the \emph{Ehrhart series} of $\mathcal P$.
Later, Stanley developed this subject, extending Ehrhart's study to rational polytopes
(see \cite{BeckRobins-book} for an exposition on the subject).

Our next goal will be to relate the spectral generating function $F_{\Gamma}(z)$ of a lens space $L:=L(q;s_1,\dots,s_n)$ with the Ehrhart series of a rational polytope naturally associated to the congruence lattice $\mathcal L:=\mathcal L(q;s_1,\dots,s_n)$.
For simplicity, we may and will assume that $s_1=1$ without loosing generality.

Let $\mu=(a_1,\dots,a_n)\in \mathcal L$.
Since $\sum_{j=1}^n s_ja_j\equiv0\pmod q$, there exists $h\in\Z$ such that $a_1=hq-a_2s_2 -\dots- a_ns_n $.
Hence
\begin{equation}
\mu = h(q,0,\dots,0)+ a_2 (-s_2,1,0,\dots,0) +\dots+ a_n(-s_n ,0\dots,0,1).
\end{equation}
We conclude that
$
\{v_1=(q,0,\dots,0),\; v_2:=(-s_2,1,0,\dots,0), \dots,\; v_n=(-s_n,0\dots,0,1)\}
$
is a $\Z$-basis of $\mathcal L$.
Let $T$ denote the matrix of change of basis to the canonical basis $\mathcal C:=\{e_1,\dots,e_n\}$, that is, $Tv_j=e_j$ for all $j$.
One clearly has $T(\mathcal L)=\Z^n$ and
\begin{equation}\label{eq:T}
T= \frac1q
\begin{pmatrix}
1& s_2  & s_3 &\dots && s_n \\
0 & q & 0&\dots && 0 \\
0&0&q& & &\vdots\\
\vdots&\vdots&&\ddots \\
0&0&&&q&0\\
0&0&0&\dots&0&q
\end{pmatrix}.
\end{equation}

We consider the polytope $\mathcal P=\{x\in \R^n: \norma{x}=\sum_{j=1}^n |x_j|\leq 1\}$ (the \emph{one-norm ball} of radius $1$ centered at the origin, also called the \emph{cross-polytope}) and $\mathcal P_{\mathcal L} =T(\mathcal P)$.
Clearly, for any $k\in\N$,
\begin{equation}
\mu\in \mathcal L\text{ and }\norma{\mu}=k
\Longleftrightarrow
\mu\in \mathcal L  \cap\big(k\mathcal P\smallsetminus (k-1)\mathcal P\big)
=  T^{-1}\big(\Z^n\cap( k\mathcal P_{\mathcal L} \smallsetminus (k-1)\mathcal P_{\mathcal L})\big) .
\end{equation}
Hence
\begin{align}
\vartheta_{\mathcal L}(z)
&= \sum_{k\geq0} N_{\mathcal L}(k)\, z^k
= \sum_{k\geq0} \#\Big( \big(k\mathcal P_{\mathcal L}\cap \Z^n\big) \smallsetminus \big((k-1)\mathcal P_{\mathcal L}\cap \Z^n\big) \Big)\, z^k \\
&= (1-z) \sum_{k\geq0} \#\left(k\mathcal P_{\mathcal L}\cap \Z^n\right)\, z^k  \notag  = (1-z)\Ehr_{\mathcal P_{\mathcal L}} (z).
\end{align}

Although $\mathcal P_{\mathcal L}$ depends on the chosen basis of the lattice $\mathcal L$, the Ehrhart series
$\Ehr_{\mathcal P_{\mathcal L}} (z)$ just depends on $\mathcal L$.
In light of \eqref{eq:F_L(z)}, we have the following result.
\begin{proposition}
Let $L=L(q;1,s_2,\dots,s_n) = \Gamma\ba S^{2n-1}$ be a lens space, and let $\Ehr_{\mathcal P_{\mathcal L}}(z)$ be the Ehrhart series of the polytope $\mathcal P_{\mathcal L}$ associated above.
Then,
\begin{equation*}
F_{\Gamma}(z) = \frac{(1-z)} {(1-z^2)^{n-1}}\; \Ehr_{\mathcal P_{\mathcal L}}(z).
\end{equation*}
\end{proposition}

Furthermore, it is well known that the Ehrhart series of an arbitrary rational polytope $\mathcal P$ can be  written as $\Ehr_{\mathcal P}(z) = (1-z^q)^{-(n+1)}\, g(z)$ for some polynomial $g(z)$ of degree less than $q(n+1)$, where $n$ is the dimension of $\mathcal P$ and $q$ is the smallest integer such that $q\mathcal P$ is integral (see for instance \cite[\S3.7]{BeckRobins-book}).

\begin{theorem}
Let $L=L(q;1,s_2,\dots,s_n) = \Gamma\ba S^{2n-1}$ be a lens space.
Then, there exists a polynomial $g_{L}(z)$ of degree less than $q(n+1)$ such that
$$
F_{\Gamma}(z) = \frac{(1-z)\; g_L(z)} {(1-z^2)^{n-1}\; (1-z^q)^{n+1}}.
$$
\end{theorem}

The previous formula implies that the $0$-spectrum of $L$ is completely determined by the polynomial $g_L(z)$, which has degree less than $q(n+1)$.
Thus, the full spectrum of $L$ is determined by finitely many numbers.

To a polytope in $\R^n$ one can canonically associate a toric variety. Hence, given a lens space $L$, one has a toric variety $X$. This variety $X$ depends on the choice of $\Z$-basis of the lattice $\mathcal L$, but different choices of $\Z$-basis give equivalent toric varieties.
As an application, Mohades and Honari found consequences of the isospectrality of lens spaces $L$ and $L'$, in terms of properties
of the associated toric varieties $X$ and $X'$,
showing in particular that isospectrality implies the equality of the dimensions of the spaces of global sections of the line bundles $\mathcal B^k$ for every $k\in\N$, where $\mathcal B$ is a suitable line bundle over $X$ (\cite[Thm.~5.3]{MohadesHonari17}).
It is now natural to investigate the validity of the converse assertion, that is, give conditions on the toric varieties implying the isospectrality of the lens spaces.

Also, we have the following interesting problem.

\begin{problem}\cite[Ques.~2]{MohadesHonari17} Extend the previous result
to lens spaces $p$-isospectral for all $p$.
\end{problem}

To conclude this section, we mention some very interesting work by Dryden, Guillemin and Sena-Dias (\cite{DGS1}, \cite{DGS2})
relating the equivariant spectrum of certain toric varieties to the geometry of the associated polytope.

\section{Related results}
In this section we give a brief summary of several recent results on spectra of lens orbifolds, not yet discussed above.

\subsection{Dirac operator}

Boldt~\cite{Boldt15} (in his Diploma Thesis \cite{Boldt-diplomarbeit}) considered the extension to the Dirac operator of Ikeda and Yamamoto's result on the non-existence of $3$-dimensional (Laplace) isospectral lens spaces.
He showed that if two $3$-dimensional lens spaces having fundamental group of prime order are Dirac isospectral, then they are isometric.
The proof follows an observation by Yamamoto in \cite{Yamamoto80} that involves $p$-adic numbers over cyclotomic fields of prime order.
It also uses an expression of the two generating functions associated to the Dirac spectrum of a spherical space form due to B\"ar~\cite[\S3]{Bar96} (see also \cite{Bar92}).

For higher dimensional spin lens spaces, Boldt and the first named author applied the method developed in \cite{LMR-onenorm} to the Dirac operator, obtaining analogous results to Theorems~\ref{thm:Fvartheta}--\ref{thm:3.6(i)} (see \cite{BoldtLauret-onenormDirac}).
In addition, several families of Dirac isospectral spin lens spaces were constructed.
Also, they obtained a finiteness condition for two spin lens spaces being Dirac-isospectral.
This allowed them to make computations, which gave strong evidences for the following claim (\cite[Conj.~6.3]{BoldtLauret-onenormDirac}, \cite[Conj.~2.14]{BoldtSchueth18})
\begin{conjecture}
	Two Dirac isospectral spin lens spaces of dimension $d\equiv 1\pmod 4$ are isometric.
\end{conjecture}

See \cite[\S2]{BoldtSchueth18} for a more detailed explanation of these results.

The assertion in the conjecture does not hold for arbitrary spin manifolds.
In fact, in~\cite{MP06}, families of pairwise Dirac isospectral compact flat manifolds in every dimension greater than or equal to $4$ are constructed.

\begin{problem}\cite[6.6]{BoldtLauret-onenormDirac}
Can the method by DeFord and Doyle~\cite{DeFordDoyle14} be applied to show Dirac isospectrality of spin lens spaces?
\end{problem}

\subsection{Computational results} \label{subsec:computationalstudy}
There have been a number of results on $p$-isospectrality of lens spaces obtained by computational methods.
Gornet and McGowan~\cite{GornetMcGowan06} used the spectral generating function in \eqref{eq:F_Gamma^p} defined by Ikeda,
to show the existence of scattered pairs of $p$-isospectral lens spaces for some $p>0$, which are not $0$-isospectral.
We refer to \cite[\S2.5]{Lauret-computationalstudy} for an analysis of these results.

In \cite{Lauret-computationalstudy}, for small values of $n$ and $q$ and for any subset $I$ of $\{0,1,\dots,n-1\}$, the computer determines the families of lens orbifolds (up to isometry) that are mutually $p$-isospectral for all $p\in I$ and such that, for each $p\notin I$, there are two elements in the family that are not $p$-isospectral.

The numerical results suggested a series of facts that were proved or conjectured.
For example,
if the lens spaces $L(q;s_1,\dots,s_n)$ and $L(q;s_1',\dots,s_n')$ are $0$-isospectral, then their corresponding covers of the same degree are also $0$-isospectral, that is, $L(q/d;s_1,\dots,s_n)$ and $L(q/d;s_1',\dots,s_n')$ are also $0$-isospectral for every  divisor $d$ of $q$ (\cite[Thm.~4.3]{Lauret-computationalstudy}).
The numerical computations seem to imply that the same fact is true for lens orbifolds, but the same proof does not extend, at least in a direct way.
\begin{conjecture} \cite[Rem.~4.4]{Lauret-computationalstudy}
	If the lens orbifolds $L(q;s_1,\dots,s_n)$ and $L(q;s_1',\dots,s_n')$ are $0$-isospectral, then $L(q/d;s_1,\dots,s_n)$ and $L(q/d;s_1',\dots,s_n')$ are also $0$-isospectral for every divisor $d$ of~$q$.
\end{conjecture}

\subsection{Spectral rigidity for 3-dimensional lens orbifolds}\label{subse:BariHunsicker}
Bari and Hunsicker~\cite{ShamsHunsicker17} considered the problem of extending Ikeda and Yamamoto's
result on $3$-dimensional lens spaces, to the orbifold case, proving that
two $3$-dimensional (and $4$-dimensional) isospectral lens orbifolds are isometric.
Thus, Conjecture~4.6 in \cite{Lauret-spec0cyclic} follows.

We next explain the steps in the argument in~\cite{ShamsHunsicker17} in the case of dimension~$3$.
Assume that $L=L(q;s_1,s_2)=\Gamma\ba S^3$ and $L'=L(q';s_1',s_2')=\Gamma'\ba S^3$ are isospectral.
We also assume that $\gcd(q,s_1,s_2)=1=\gcd(q',s_1',s_2')$.
A simple observation shows that $q=q'$.
The proof is naturally divided into the following five cases:
\begin{enumerate}
\item $L$ and $L'$ are lens spaces;
\item exactly one of $L$ and $L'$ is a lens space;
\item $\gcd(q,s_1)=1=\gcd(q,s_1')$ and $\gcd(q,s_2)\neq 1\neq \gcd(q,s_2')$;
\item $\gcd(q,s_1)=1$, and $\gcd(q,s_2)$, $\gcd(q,s_1')$, $\gcd(q,s_2')$ are greater than $1$;
\item $\gcd(q,s_j)$ and $\gcd(q,s'_j)$ are greater than $1$ for $j=1,2$.
\end{enumerate}
We give the main ideas in each case.
Case (1) was proved in~\cite{Yamamoto80}.
The second case follows from the following general fact proved in~\cite[Prop.~3.4(ii)]{GordonRossetti03}: if two isospectral orbifolds share the same
Riemannian universal cover, and one of them is a manifold, then the other is also a manifold.

In the remaining cases, the authors used the spectral invariants given by the residues of the spectral generating function $F_{\Gamma}(z)$ at the $q$-th primitive roots of unity, as done in~\cite{IkedaYamamoto79} and~\cite{Yamamoto80}.
For example, in case (3), after  assuming $s_1=s_1'=1$, they proved that the residues of $F_{\Gamma}(z)$ and $F_{\Gamma'}(z)$ at $\xi_q=e^{2\pi \mi /q}$ are equal to
\begin{equation}
\frac{-2\, \xi_q}{q(1-\xi_q^{1-s_2}) (1-\xi_q^{1+s_2})}
=
\frac{-2\, \xi_q}{q(1-\xi_q^{1-s_2'}) (1-\xi_q^{1+s_2'})} .
\end{equation}
By straightforward cancellations, one obtains that $\xi_q^{s_2}+\xi_q^{-s_2} = \xi_q^{s_2'}+\xi_q^{-s_2'}$, thus $\cos(2\pi s_2/q)=\cos(2\pi s_2'/q)$ and therefore $s_2\equiv  \pm s_2'\pmod q$.
The last condition implies that $L(q;1,s_2)$ and $L(q;1,s_2')$ are isometric.

Case (4) is even simpler since $\xi_q$ is a pole of $F_{\Gamma}(z)$ but not of $F_{\Gamma'}(z)$.
Case (5) is the challenging one, and the proof in \cite{ShamsHunsicker17} is rather difficult to check.

\subsection{Quotient of other compact symmetric spaces} \label{subsec:spec0cyclic}
The general approach introduced in Section~\ref{sec:approach} works also for other normal homogeneous spaces besides odd-dimensional spheres.
We let $(G,K)$ be such that $G/K$ is a compact symmetric space of real rank one.
 As in the case of the sphere, the corresponding spherical representations (i.e. those in  $\widehat G_K$) are given by a family parametrized by non-negative integers, say $\widehat G_{\tau_0} = \{(\sigma_k,V_{\sigma_k})\}_{k\geq0}$ (see \cite[\S3.1]{Lauret-spec0cyclic} for further details), with associated eigenvalues satisfying $\lambda(C,\sigma_k)<\lambda(C,\sigma_{k+1})$ for all $k\geq0$. 
In particular, every eigenspace is an irreducible $G$-module.

We assume that $\Gamma$ is a finite subgroup of the maximal torus $T$ of $G$.
Proposition~\ref{prop:tau-spectrum} tells us that every  eigenvalue of the Laplace--Beltrami operator on $\Gamma\ba G/K$ is of the form $\lambda(C,\sigma_k)$ for some $k\geq0$ with multiplicity equal to $\dim V_{\sigma_k}^\Gamma \dim \Hom_K(V_{\sigma_k},W_{\tau_0})$.
It remains to determine $\dim V_{\sigma_k}^\Gamma$ for every $k$, since the number $\dim \Hom_K(V_{\sigma_k},W_{\tau_0})$ is known in this case (it is equal to one for $G/K\simeq S^{2n-1}, S^{2n}, P^n(\C)$).
We define an analogue of Ikeda's spectral generating function \eqref{eq:spectralgeneratingfunction} by
\begin{equation}
F_\Gamma(z)=F_{G,K,\Gamma}(z) = \sum_{k\geq0} \dim {V_{\sigma_k}^\Gamma}\; z^k.
\end{equation}
From \eqref{eq:dimV_pi^Gamma}, we have that
$\dim V_{\sigma_k}^\Gamma =\sum_{\mu \in \mathcal L_\Gamma} m_{\sigma_k}(\mu)$, where
$
\mathcal L_\Gamma$ is the lattice defined in \eqref{eq:mathcalL_Gamma}.

Usually, it is not difficult to describe $\mathcal L_\Gamma$ in terms of the parameters of the generators of $\Gamma$.
The main difficulty is to obtain an expression for the weight multiplicity $m_{\sigma_k}(\mu)$ depending only on $\|\mu\|$ for an adequate norm $\|\cdot\|$ on the corresponding weight lattice $P(G)$.

In \cite[Prop.~3.2]{Lauret-spec0cyclic}, such an expression is given for $(G,K)$ equal to
$(\SO(2n),\SO(2n-1))$,
$(\SO(2n+1),\SO(2n))$,
$(\SU(n+1),\mathrm{S}(\Ut(n)\times\Ut(1)))$
and $(\Sp(2),\Sp(1)\times\Sp(1))$.
As a consequence, one obtains an expression analogous to  \eqref{eq:dimHH_k^Gamma3} for the multiplicity $\dim V_{\sigma_k}^\Gamma$ of the eigenvalue $\lambda(C,\sigma_k)$ of the Laplace--Beltrami operator on $\Gamma\ba G/K$ (see \cite[Thm.~3.4]{Lauret-spec0cyclic}), as well as an expression for $F_{\Gamma}(z)$.

For instance, if $G/K\simeq S^{2n}$, the multiplicity of the eigenvalue $\lambda(C,\sigma_k)=k(k+2n-1)$ is given by
\begin{equation}\label{eq:dimV_sigma_k}
\dim V_{\sigma_k} ^\Gamma
= \sum_{r=0}^{\lfloor k/2\rfloor }  \tbinom{r+n-1}{n-1} \left( N_{\mathcal L_\Gamma}(k-2r)+N_{\mathcal L_\Gamma}(k-1-2r) \right)
\end{equation}
and
\begin{equation}\label{eq:F_GammaSO(2n+1)SO(2n)}
F_{\Gamma}(z) = \frac{1+z}{(1-z^2)^n} \vartheta_{\mathcal L_{\Gamma}}(z).
\end{equation}
Here, since the standard maximal tori $T$ of $\SO(2n+1)$ and $\SO(2n)$ coincide, as well as their corresponding weight lattices (which are identified with $\Z^n$), the norm $\norma{\cdot}$ and the function $N_{\mathcal L_{\Gamma}}(k)$ are defined as in Subsection~\ref{subsec:alternative}

Laplace operators on $\Gamma\ba G/K$ twisted by a character $\chi$ of $\Gamma$ were also considered in \cite{Lauret-spec0cyclic}.
In this case, the associated `lattice' $\mathcal L_{\Gamma,\chi}$ turns out to be a shifted lattice of the original congruence lattice $\mathcal L_\Gamma$, that is, $\mathcal L_{\Gamma,\chi} = \mu_0+\mathcal L_{\Gamma}$  for some $\mu_0\in P(G)$.
The remaining formulas (e.g.\ \eqref{eq:dimHH_k^Gamma3}, \eqref{eq:F_L(z)2}, \eqref{eq:dimV_sigma_k} \eqref{eq:F_GammaSO(2n+1)SO(2n)}) hold in this generality, by replacing $\mathcal L_\Gamma$ by $\mathcal L_{\Gamma,\chi}$.

In \cite[\S4]{Lauret-spec0cyclic} many tables are given with computational calculations providing many isospectral examples with respect to twisted Laplace operators on lens orbifolds.
For instance, it was found that for some characters,  the $3$-dimensional lens orbifolds $L(4;1,2)$ and $L(4;0,1)$ are isospectral with respect to the corresponding twisted Laplace operators, contrasting with the non-existence in the untwisted case in dimension $3$.
Another peculiar example is that the lens space $L(5;1,2)$ with any non-trivial twist is isospectral to the lens orbifold $L(5;0,1)$ with a certain non-trivial twist.
Examples of twisted isospectrality between a manifold and an orbifold (with singularities) were already constructed by the last two named authors in \cite{MR-comparison}, where the case of compact flat manifolds was studied.

These examples confirm that there is not so much geometric information encoded in the spectral information associated to the Laplace operator twisted by a (non-trivial) character of the fundamental group.

\subsection{Harmonic-counting measure}
Let $\mathcal L=\mathcal L(q;s_1,\dots,s_n)$ be the congruence lattice associated to the lens space $L(q;s_1,\dots,s_n)$ given by elements $(a_1,\dots,a_n)\in\Z^n$ satisfying that $\sum_{i=1}^n a_is_i$ is divisible by $q$.
Mohades and Honari in \cite{MohadesHonari16} defined the \emph{harmonic-counting measure associated to $\mathcal L$} as the measure $\nu_{\mathcal L}$ on the Borel $\sigma$-algebra of the unit sphere $S^{n-1}$ as
\begin{equation}
\nu_{\mathcal L}(U) = \lim_{t \to\infty}
\frac{
	\displaystyle
	\sum_{k=0}^{\lfloor t\rfloor} \sum_{r=0}^{\lfloor k/2 \rfloor } \binom{r+n-2}{n-2} \,\#\{\mu \in \mathcal L\cap C(U): \norma{\mu}=k-2r\}
}{t^{2n-1}},
\end{equation}
where $C(U)= \{a\mu: \mu\in U,\, a>0\}$ (the cone in $\R^n$ induced by $U$).
When $U=S^{n-1}$, \eqref{eq:dimHH_k^Gamma} tells us that the numerator on the right hand side in the above formula equals the number of eigenvalues (counted with multiplicities) of the Laplace--Beltrami operator on $L$ less than or equal to $\lambda_{\lfloor t\rfloor}=\lfloor t\rfloor(\lfloor t\rfloor+2n-2)$.

By using the Weyl law for the lens space $L(q;s_1,\dots,s_n)$, they proved that $\nu_{\mathcal L}$ is a finite measure with total value  $\nu_{\mathcal L}(S^{n-1}) = \frac{2}{q(2n-1)!}$ and such that
\begin{equation}
\nu_{\mathcal L}(U)= \frac{  1}{(n-2)!2^{n-1}} \op{Vol}(T(C(U))\cap \mathcal B)
\int_{0}^1 x^{n-2} (1-x)^{n} dx,
\end{equation}
where $T$ is the generating matrix of $\mathcal L$ defined in \eqref{eq:T}  and $\mathcal B$ is the unit one-norm ball.

\bibliographystyle{plain}

\end{document}